\documentclass[reqno,11pt,a4paper]{amsart}

\usepackage{leftidx}
\usepackage{xcolor}
\usepackage{amssymb,amsmath}
\usepackage{enumerate}

\newtheorem{MainTheorem}{Theorem}
\newtheorem{Proposition}{Proposition}[section]

\newtheorem{Lemma}[Proposition]{Lemma}
\newtheorem{Theorem}[Proposition]{Theorem}
\newtheorem{Corollary}[Proposition]{Corollary}

\DeclareMathOperator{\Val}{Val}

\DeclareMathOperator{\vol}{vol}

\DeclareMathOperator{\Kl}{Kl}
\DeclareMathOperator{\Gr}{Gr}
\DeclareMathOperator{\tr}{tr}
\DeclareMathOperator{\Id}{Id}

\DeclareMathOperator{\im}{Im}
\DeclareMathOperator{\re}{Re}
\DeclareMathOperator{\diag}{diag}

\DeclareMathOperator{\rank}{rank}
\DeclareMathOperator{\Sp}{Sp}
\DeclareMathOperator{\SO}{SO}
\DeclareMathOperator{\GL}{GL}
\DeclareMathOperator{\SU}{SU}

 \newcommand{\spsp}{\Sp(2)\Sp(1)}

\newcommand{\spann}{\mathrm{span}}
\newcommand{\R}{\mathbb{R}}

\newcommand{\h}{\mathbb{H}}
\renewcommand{\i}{\mathbf{i}}
\renewcommand{\j}{\mathbf{j}}
\renewcommand{\k}{\mathbf{k}}
\title[Valuations on the quaternionic plane]{Classification of invariant valuations on the quaternionic plane} 
\author{Andreas Bernig}
\author{Gil Solanes}
 
\email{bernig@math.uni-frankfurt.de}
\email{solanes@mat.uab.cat}

\address{Institut f\"ur Mathematik, Goethe-Universit\"at Frankfurt,
Robert-Mayer-Str. 10, 60054 Frankfurt, Germany}
\address{Departament de Matem\`atiques, Universitat Aut\`onoma de Barcelona, 08193 Bellaterra, Spain}

\thanks{A.B. was supported by DFG grants BE 2484/3-1 and BE 2484/5-1. G.S. was supported by FEDER/MEC grant MTM2009/07594 and AGAUR grant SGR2009-1207.\\ AMS 2010 {\it Mathematics subject
classification}: 53C65, 
53C35. 
}

\begin{document}

\begin{abstract}
We describe the orbit space of the action of the group $\spsp$ on the real Grassmann manifolds
$\Gr_k(\h^2)$ in terms of certain quaternionic matrices of Moore rank not larger than $2$. We
then 
give a complete classification of valuations on the quaternionic plane $\h^2$ which are invariant under the action
of the group $\spsp$. 
\end{abstract}

\maketitle

\section{Introduction and statement of results}
\subsection{Background}
A valuation is a finitely additive map from the space of compact convex subsets of some vector space into an abelian
semi-group. Since Hadwiger's famous characterization of (real-valued) continuous valuations which are euclidean motion
invariant, classification results for valuations have long played a prominent role in convex and integral
geometry. 

Many generalizations of Hadwiger's theorem were obtained recently. On the one hand, valuations with values in
some abelian semi-group other than the reals were characterized. The most important examples are tensor valuations
\cite{alesker99, hug_schneider_schuster_a, hug_schneider_schuster_b, mcmullen97}, Minkowski valuations
\cite{abardia12, abardia_bernig, haberl10, ludwig_2005, schuster10, schuster_wannerer}, curvature measures
\cite{bernig_fu_solanes,schneider78} and area measures \cite{wannerer_area_measures, wannerer_unitary_module}.  On the
other hand, invariance with respect to the
euclidean group was weakened to invariance with respect to translations or rotations only \cite{alesker99_annals,
alesker_mcullenconj01}, or with respect to a
smaller group of isometries. Next we briefly describe the main results in this line.

Let $V$ be a finite-dimensional vector space and $G$ a group
acting linearly on $V$. The space of {scalar-valued}, $G$-invariant,
translation invariant continuous valuations on $V$ will be denoted by $\Val^G$. Hadwiger's theorem applies in
the case where $V$
is a euclidean vector space of dimension $n$, and $G=\SO(V)$. It states that $\Val^G$ is spanned by the so-called
intrinsic volumes $\mu_0,\ldots,\mu_n$. In particular, $\Val^{\SO(V)}$ is finite-dimensional. From this fact, one can
easily derive integral-geometric formulas like Crofton formulas and kinematic formulas \cite{klain_rota}. 

In the same spirit, kinematic formulas with respect to a smaller group $G$ exist provided that $\Val^G$ is
finite-dimensional. Although it is known which groups have this property, much less is known about the explicit form of
such formulas. Alesker \cite{alesker_survey07} has shown that $\Val^G$ is finite-dimensional if and only if $G$ acts
transitively on the unit sphere. Such groups were classified by Montgomery-Samelson \cite{montgomery_samelson43} and
Borel \cite{borel49}. There are six infinite lists
\begin{equation} \label{eq_transitive_groups}
 \mathrm{SO}(n), \mathrm{U}(n), \mathrm{SU}(n), \mathrm{Sp}(n), {\mathrm{Sp}(n)\mathrm{U}(1), \mathrm{Sp}(n)
\mathrm{Sp}(1)}
\end{equation}
and three exceptional groups
\begin{equation} \label{eq_exceptional_groups}
 \mathrm{G}_2, \mathrm{Spin}(7), \mathrm{Spin}(9). 
\end{equation}

The euclidean case is $G=\mathrm{SO}(n)$ where Hadwiger's theorem applies. In the hermitian case $G=\mathrm{U}(n)$
or $G=\mathrm{SU}(n)$, recent results have revealed a lot of unexpected algebraic structures yielding a relatively
complete picture \cite{abardia_gallego_solanes, alesker_mcullenconj01, bernig_fu_hig,
bernig_fu_solanes, park02, tasaki03}. Hadwiger-type theorems for the groups $\mathrm{G}_2$ and $\mathrm{Spin}(7)$ are
also known \cite{bernig_g2}. In the
remaining cases, i.e. the quaternionic cases $G=\mathrm{Sp}(n)$, $G=\mathrm{Sp}(n)\mathrm{U}(1)$ and
$G=\mathrm{Sp}(n)\mathrm{Sp}(1)$ as well as in the case $G=\mathrm{Spin}(9)$, only the dimension of $\Val^G$ is
known \cite{bernig_qig, voide}. 

The combinatorial formulas from \cite{bernig_qig} indicate that the classification of invariant valuations on
quaternionic vector spaces will be a rather subtle subject. Note that the case $n=1$ can be reduced to the hermitian
case, since $\Sp(1)=\mathrm{SU}(2)$. For higher dimensions, not much is known, except the construction of one example of
an $\Sp(n)\Sp(1)$-invariant valuation by Alesker \cite{alesker05}. 

\subsection{Results of the present paper}
In this article, we establish a complete Hadwiger-type theorem for the group $\Sp(2)\Sp(1)$ acting on the
two-dimensional quaternionic space $\h^2$. More precisely, we find an explicit basis of the space of invariant
valuations $\Val^{\Sp(2)\Sp(1)}$.  The  description of the basis is given in terms of Klain functions, which are
invariant functions on the real Grassmannians of $\h^2$. 

Our first main theorem concerns the
orbit space of the action of $\Sp(2)\Sp(1)$ on the real Grassmann manifolds
$\Gr_k:=\Gr_k(\h^2)$. It is
formulated in terms of the Moore rank of hyperhermitian matrices, whose definition will be recalled in the
next section.
Since taking orthogonal complements commutes with the action of $\Sp(2)\Sp(1)$, it will be enough to consider the
case $k
\leq 4$.  

\begin{MainTheorem} \label{thm_orbits}
 Let $2 \leq k \leq 4$. Given a tuple of real numbers $\lambda_{pq}, 1 \leq p < q \leq k$ we define the quaternionic
hermitian matrix $M_\lambda$ by 
\begin{displaymath}
 M_\lambda:=\begin{cases} \left(\begin{array}{c c} 1 & \lambda_{12} \i \\
             -\lambda_{12} \i & 1 
            \end{array}\right) & k =2\\
\left(\begin{array}{c c c} 1 & \lambda_{12} \i & \lambda_{13} \j \\
             -\lambda_{12} \i & 1 & {-\lambda_{23}} \k \\
-\lambda_{13}\j & \lambda_{23} \k & 1 
            \end{array}\right) & k=3\\
\left(\begin{array}{c c c c} 1 & \lambda_{12} \i & \lambda_{13} \j & \lambda_{14} \k \\
             -\lambda_{12} \i & 1 & {-\lambda_{23}} \k & \lambda_{24} \j\\
-\lambda_{13} \j & \lambda_{23} \k & 1 & -\lambda_{34} \i\\
-\lambda_{14} \k & -\lambda_{24} \j & \lambda_{34} \i & 1
            \end{array}\right) & k=4.
\end{cases}
\end{displaymath}
Let {$\mathbb{Z}_2^k$ and the permutation group} $\mathcal{S}_k$ act on such a tuple by 
\begin{align}
  (\epsilon\cdot\lambda)_{p,q} & :=\epsilon_p \epsilon_q \lambda_{pq},&\epsilon\in \mathbb{Z}_2^k\label{actionZ}\\
 (\sigma\cdot\lambda)_{p,q} & := \lambda_{\sigma(p)\sigma(q)}=\lambda_{\sigma(q)\sigma(p)},&\sigma\in \mathcal{S}_k.\label{actionS}
\end{align}
Then the quotient
$\Gr_k/\Sp(2)\Sp(1)$ is of dimension $(k-1)$ and homeomorphic to the quotient 
\begin{displaymath}
 X_k:=\{\lambda_{pq} \in [-1,1], 1 \leq p<q \leq k: \rank M_\lambda \leq 2\} / \mathbb{Z}_2^k \times
\mathcal{S}_k. 
\end{displaymath}
The orbit corresponding to {$[\lambda] \in X_k$} contains a plane $V$ admitting a basis $v_1,\ldots,v_k$ such
that
\begin{displaymath}
K(v_i,v_j)=(M_\lambda)_{i,j}\qquad i,j=1,\ldots,k,
\end{displaymath}
where $K$ is the quaternionic hermitian product of $\h^2$.
\end{MainTheorem}

The construction of this homeomorphism is roughly as follows. Given a plane $V \in \Gr_k(\h^2)$, we construct an orthonormal basis $v_1,\ldots,v_k$ of $V$ such that the matrix $Q=(K(v_i,v_j))$ has a special shape: if $k \in \{3,4\}$, the pure quaternions $q_{12},q_{13},q_{23}$ are pairwise orthogonal, and moreover  $q_{12}\|q_{34}, q_{13}\|q_{24}, q_{14}\|q_{23}$ if $k=4$. Then $Q$ is $\Sp(1)$-conjugate to a matrix $M_\lambda$, and $V$ is mapped to $[\lambda]$.

Note that the condition on the Moore rank is a system of polynomial equations in the $\lambda_{pq}$, which can be written down explicitly using  equations \eqref{eq_moore3} and \eqref{eq_moore4}. 

\begin{Corollary} \label{cor_param_angles}
 Every $\Sp(2)\Sp(1)$-orbit in $\Gr_k$ contains a {$k$-plane} of the form
\begin{align*}
&\spann\{(\cos\theta_1,\sin\theta_1),(\cos\theta_2,\sin\theta_2)\mathbf i\} & k=2,\\
&\spann\{(\cos\theta_1,\sin\theta_1),(\cos\theta_2,\sin\theta_2)\mathbf i,(\cos\theta_3,\sin\theta_3)\mathbf j\} &
k=3,\\
&\spann\{(\cos\theta_1,\sin\theta_1),(\cos\theta_2,\sin\theta_2)\mathbf i,(\cos\theta_3,\sin\theta_3)\mathbf
j,(\cos\theta_4,\sin\theta_4)\mathbf k\} & k=4,
\end{align*}
where $\theta_1,\ldots,\theta_4 \in [0,2\pi]$.
The corresponding $[\lambda]\in X_k$ is given by $\lambda_{pq}=\cos(\theta_p-\theta_q)$.
\end{Corollary}

Let us now describe the Hadwiger-type theorem, which is our second main result. 
The space of continuous, translation invariant valuations on an $n$-dimensional vector space $V$ is denoted by $\Val(V)$
or just $\Val$ if there is no risk of confusion. A valuation $\phi \in \Val$ is called {\it even} if $\phi(-
B)=\phi(B)$
and {\it odd} if $\phi(-B)=-\phi(B)$ for each convex body $B$. If $\phi(tB)=t^k \phi(B)$ for all $t>0$ and all
$B$,
then $\phi$ is said to be {\it homogeneous of degree $k$}. The space of even/odd valuations of degree $k$ is denoted by
$\Val_k^\pm$. A
fundamental result by McMullen \cite{mcmullen77} is the decomposition 
\begin{displaymath}
 \Val = \bigoplus_{\substack{k=0,\ldots,n\\ \epsilon=\pm}} \Val_k^\epsilon. 
\end{displaymath}

An even, continuous and translation
invariant valuation can be described by its Klain function, which is defined as follows. Let $\phi \in \Val_k^+$ and $E
\in \Gr_k(V)$, the Grassmann manifold of $k$-planes in $V$. Then the restriction of $\phi$ to $E$ is a multiple of the
Lebesgue measure, and the corresponding factor is denoted by $\Kl_\phi(E)$. The function $\Kl_\phi \in
C(\Gr_k(V))$ is 
called the Klain function of $\phi$. The map $\Kl: \Val_k^+ \to C(\Gr_k(V))$ is in fact injective, as was shown
by Klain
\cite{klain00}.

Let us now specialize to the group $\spsp$ acting on $V=\h^2$. The dimension of the space of
$k$-homogeneous
$\spsp$-invariant valuations was computed in \cite{bernig_qig}:
\begin{equation}\label{table_dims}
\begin{array}{c | c |c | c | c | c | c | c | c | c} 
k & 0 & 1 & 2 & 3 & 4 & 5 & 6 & 7 & 8\\ \hline
\dim \Val_k^{\spsp} &  1 & 1 & 2 & 3 & 5 & 3 & 2 & 1 & 1
\end{array}
\end{equation} 
Since the group $\spsp$ contains $-\Id$, invariant valuations are even. We will characterize them in terms of
their Klain functions. To do so, consider the following invariant functions on $\Gr_k, 0 \leq k \leq 4$,
which
are defined in terms of the
coordinates $\lambda=(\lambda_{ij})$ of $\Gr_k/{\spsp}$ from Theorem \ref{thm_orbits}.  
\begin{align*}
f_{k,0}(\lambda) & := 1, \quad k=0,\ldots,4\\
f_{2,1}(\lambda) & :=\lambda_{12}^2\\
f_{3,1}(\lambda) & :=\lambda_{12}^2+\lambda_{13}^2+\lambda_{23}^2\\ 
f_{3,2}(\lambda) & :=\lambda_{12}^2\lambda_{23}^2+\lambda_{13}^2\lambda_{23}^2+ \lambda_{12}^2 \lambda_{13}^2\\
f_{4,1}(\lambda) &
:=\lambda_{12}^2+\lambda_{13}^2+\lambda_{14}^2+\lambda_{23}^2+\lambda_{24}^2+\lambda_{34}^2\\
f_{4,2}(\lambda) & :=\lambda_{12}^2\lambda_{34}^2+\lambda_{13}^2\lambda_{24}^2+\lambda_{14}^2\lambda_{23}^2\\
f_{4,3}(\lambda) &
:= \lambda_{12}^2\lambda_{13}^2+\lambda_{12}^2\lambda_{14}^2+\lambda_{13}^2\lambda_{14}^2+\lambda_{12}^2\lambda_{
23 } ^2+\lambda_{12}^2\lambda_{24}^2+\lambda_{23}^2\lambda_{24}^2\\
& \quad +
\lambda_{13}^2\lambda_{23}^2+\lambda_{13}^2\lambda_{34}^2+\lambda_{23}^2\lambda_{34}^2+\lambda_{14}^2\lambda_{24}
^2+\lambda_{14}^2\lambda_{34}^2+\lambda_{24}^2\lambda_{34}^2\\
f_{4,4}(\lambda) & :=2\lambda_{12}\lambda_{13}\lambda_{23}^2\lambda_{24}\lambda_{34}
+2\lambda_{12}\lambda_{13}\lambda_{14}^2\lambda_{24}\lambda_{34}
+2\lambda_{12}\lambda_{23}\lambda_{13}^2\lambda_{14}\lambda_{34}\\
& \quad +2\lambda_{12}\lambda_{23}\lambda_{24}^2\lambda_{14}\lambda_{34}
+2\lambda_{24}\lambda_{23}\lambda_{12}^2\lambda_{14}\lambda_{13}
+2\lambda_{24}\lambda_{23}\lambda_{34}^2\lambda_{14}\lambda_{13}\\
& \quad
+3( \lambda_{12}^2\lambda_{13}^2\lambda_{14}^2+\lambda_{12}^2\lambda_{23}^2\lambda_{24}^2+\lambda_{13}^2\lambda_{
23} ^2\lambda_{34}^2+\lambda_{14}^2\lambda_{ 24 }^2\lambda_{34}^2).
\end{align*} 

Noting that $\Gr_k \cong \Gr_{8-k}$ for all $k$, we define $f_{k,i}:=f_{8-k,i}$ for $5 \leq k \leq 8$.

\begin{MainTheorem} \label{thm_hadwiger_type}
For each $0 \leq k \leq 8$ and each $0 \leq i \leq \dim \Val_k^{\spsp}-1$, there exists a unique valuation $\phi
\in \Val^{\spsp}_k$ whose Klain function
is $f_{k,i}$. These valuations form a basis of $\Val_k^{\spsp}$.  
\end{MainTheorem}

Moreover, we will find Crofton measures for these valuations. In the proof of this theorem, we will first use differential geometric methods to show that certain linear
combinations of the functions $f_{k,i}$ are eigenfunctions of the Laplace-Beltrami operator on $\Gr_k$. Then we
will use
re\-pre\-sen\-tation-theoretic tools, in particular the recent computation of the multipliers of the $\alpha$-cosine transform
by \'Olafsson-Pasquale \cite{olafsson_pasquale}, in order to construct valuations with the given Klain functions. As a corollary to their theorem, we prove a formula for the multipliers of the classical cosine transform which might be of independent interest. To
see that the so-constructed valuations form a basis, we use the recent computation of $\dim \Val^{\spsp}$ in
\cite{bernig_qig}. 

Let us mention that Alesker \cite{alesker05} has constructed a quaternionic version of Kazarnovskii's
pseudo-volume (compare \cite{alesker03_un,kazarnovskii} for Kazarnovskii's pseudo-volume on $\mathbb{C}^n$). Given any
$n$, Alesker's pseudo-volume is a
continuous, translation invariant, $\Sp(n)\Sp(1)$-invariant valuation of degree $n$ on $\h^n$. It has the property that
its restriction to each quaternionic hyperplane vanishes. In the present case $n=2$, a quaternionic line inside $\h^2$
is given by the angles $\theta_1=\theta_2=0$, i.e. $\lambda_{12}=1$. It follows that Alesker's
pseudo-volume is a real multiple of the degree $2$ valuation with Klain function $f_{2,0}-f_{2,1}$.

\subsection*{Acknowledgments}

We would like to thank Semyon Alesker, Joseph Fu and Franz Schuster for fruitful discussions and useful remarks.

\section{Quaternionic linear algebra}

The quaternionic skew field $\h$ is defined as the real algebra generated by $1,\i,\j,\k$ with the relations
$\i^2=\j^2=\k^2=-1, \i\j\k=-1$. The conjugate of a quaternion $q:=a+b\,\i+c\,\j+d\,\k$ is defined by $\bar q:=a-b\,\i-c\,\j-d\,\k$,
its norm by $\sqrt{q\bar q}$. The quaternions of norm $1$ form the Lie group $\Sp(1)$ which is isomorphic to
$\SU(2)$.
Conjugation by an element $\xi \in \Sp(1)$ fixes the real line pointwise and acts as a rotation on the pure imaginary
part $\im \h=\R^3$, moreover all rotations are obtained in this way.

Let $V$ be a quaternionic (right) vector space of dimension $n$. We endow $V$
with a quaternionic hermitian form $K$, i.e. an $\R$-bilinear form 
\begin{displaymath}
 K:V \times V \to \h
\end{displaymath}
such that  
\begin{enumerate}[i)]
\item $K$ is conjugate $\h$-linear in the first and $\h$-linear in the second
factor, i.e. 
\begin{displaymath}
 K(v q, w r)=\bar q K(v,w)r, \quad q,r \in \h,
\end{displaymath}
\item $K$ is hermitian in the sense that 
\begin{displaymath}
 K(w,v)=\overline{K(v,w)},
\end{displaymath}
\item $K$ is positive definite, i.e. 
\begin{displaymath}
 K(v,v)>0 \quad \forall v \neq 0. 
\end{displaymath}
\end{enumerate}

The standard example of such a form is given in $V=\h^n$ by
\begin{displaymath}
K(v,w)=\sum_{i=1}^n \bar v_i w_i, \quad v=(v_1,\ldots,v_n),
w=(w_1,\ldots,w_n) \in \h^n.
\end{displaymath}
 
The group $\GL(V,\h)=\GL(n,\h)$ is defined as the group of all $\h$-linear automorphisms of $V$. The subgroup of
$\GL(V,\h)$ of all elements preserving $K$ is called the {\it compact
symplectic group} and denoted by $\Sp(V,K)$ or $\Sp(n)$. It acts from the left on
$V$. An important fact is that this action is transitive on the unit sphere in $V$. In the case $V=\h^n$, the group
$\Sp(n)$ consists of all quaternionic matrices $A$ such that $A^*A=\Id$. Here $A^*$ denotes the conjugate transpose of
$A$.

The action of $\Sp(n) \times \Sp(1)$ by left and right multiplication on $V$ has kernel
$\mathbb{Z}_2=\{(\Id,1),(-\Id,-1)\}$. The quotient group
is denoted by $\Sp(n)\Sp(1)$. It acts effectively on $V$.  
  
Let $Q=(q_{ij})$ be a quaternionic $n \times n$ matrix. Viewing $\h^n$ as a {\it right} $\h$-vector space, $Q$
acts as a quaternionic linear map $Q:\h^n \to \h^n$ by multiplication from the left. Writing $\h=\R^4$, we obtain
a corresponding real linear
map $\leftidx{^\R}{Q}:\R^{4n} \to \R^{4n}$.  

A square matrix $Q$ with quaternionic entries is called {\it hyperhermitian} if $Q^*=Q$, i.e. $q_{ji}=\bar q_{ij}$ for
all $i,j$. In particular, the diagonal entries are real. The determinant of $\leftidx{^\R}{Q}$ is a polynomial of
degree $4n$ in the
$n(2n-1)$ real components of $Q$. The {\it Moore determinant} is the unique polynomial {$\det$} of degree $n$ in
the same
variables which satisfies $\det(Q)^4=\det(\leftidx{^\R}{Q})$ and $\det(\Id)=1$. Note that the Moore determinant
is defined
only on hyperhermitian matrices. We refer to \cite{alesker03, alesker05, aslaksen96} for more information on the Moore
determinant and its relation to other determinants of quaternionic matrices such as the Dieudonn\'e determinant. 

If $Q$ is a hyperhermitian matrix, there exists a matrix $A \in \Sp(n)$ and a
diagonal
matrix $D$ with real entries such that $Q=A^* D A$. Then $\det(Q)=\det(D)$. The diagonal entries in $D$ are the
(Moore-) eigenvalues of $Q$. {More generally, if $Q$ is hyperhermitian and $A$ is any quaternionic matrix, then 
\begin{displaymath}
 \det(A^*QA)=\det Q \det(A^*A),
\end{displaymath}
compare \cite{alesker05}, Thm. 1.2.9.} 

The {\it Moore rank} of $Q$ is the quaternionic dimension of the image of $Q$, or equivalently the number of non-zero
eigenvalues. Clearly the Moore rank is maximal if and only if $\det(Q) \neq 0$. 

We will need explicit formulas for Moore determinants of small size which can be computed using the results from
\cite{aslaksen96}. For $M_\lambda$ as in Theorem \ref{thm_orbits}, the Moore determinant is given by
\begin{align}
\det M_\lambda &=1-\lambda_{12}^2, \hspace{7cm}k=2,\\
\label{eq_moore3} \det M_\lambda&=1-\lambda_{12}^2-\lambda_{13}^2-\lambda_{23}^2+2\lambda_{12} \lambda_{13}
\lambda_{23},\hspace{2.7cm} k=3,\\
\label{eq_moore4} \det
M_\lambda&=1-\lambda_{12}^2-\lambda_{13}^2-\lambda_{14}^2-\lambda_{23}^2-\lambda_{24}^2-\lambda_{34}^2 \nonumber\\
& \quad +2\lambda_{23} \lambda_{34}
\lambda_{24}+2\lambda_{12}\lambda_{23}\lambda_{13}+2\lambda_{12}
\lambda_{24}\lambda_{14} +2\lambda_{13}\lambda_{34}\lambda_{14} \\
& \quad +\lambda_{12}^2\lambda_{34}^2+\lambda_{23}^2\lambda_{14}^2+\lambda_{13}^2\lambda_{24}^2 \nonumber\\
& \quad -2\lambda_{12}\lambda_{23}\lambda_{
34}\lambda_{14}-2\lambda_{12}\lambda_{24}\lambda_{13}\lambda_{34}-2\lambda_{13}
\lambda_{24}\lambda_{23}\lambda_{14},&k=4.\nonumber
\end{align}
For $k=4$, the Moore determinants of the diagonal $3\times 3$ submatrices of $M_\lambda$ can be computed by
\eqref{eq_moore3} since $\det$ is invariant under $\Sp(1)$-conjugation.

\section{Grassmann orbits}\label{section_orbits}

The aim of this section is the description of the orbit spaces of the action of the group {$G:=\Sp(2)\Sp(1)$} on the
Grassmann spaces
$\Gr_k$. Note that $\Gr_k \cong \Gr_{8-k}$, so we may assume $k \leq 4$. In the
cases $k=0,1$, the
action is transitive, so we are left with  $k=2,3,4$. Theorem \ref{thm_orbits} will follow from Theorems \ref{thm_2_orbit}, \ref{thm_3_orbit} and \ref{thm_4_orbit} below.

The following propositions will be useful.
\begin{Proposition}\label{prop_surj}
Let $Q=(q_{ij})$ be a $k\times k$ hyperhermitian matrix with Moore rank at most $2$ and non-negative
eigenvalues. Then there exist $u_{1},\ldots,u_{k}\in\h^2$ such that
\begin{displaymath}
K(u_{i},u_{j})=q_{ij}\qquad {\forall i,j.}
\end{displaymath}
\end{Proposition}

\proof
We may decompose $Q=A^* D A$ where $A=(a_{ij})\in \Sp(k)$ and 
$D= \diag(\delta_{1},\delta_{2},0,\ldots,0)$. Then
\begin{displaymath}
u_i=\left(\sqrt{\delta_{1}}a_{1i},\sqrt{\delta_{2}}a_{2i}\right) \in \h^2,\qquad i=1,\ldots,k
\end{displaymath}
are such that $K(u_i,u_{j})=q_{ij}$ for all $i,j$.
\endproof

\begin{Proposition}\label{prop_inj}
Let $u_1,\ldots,u_k \in \h^n$ and $v_1,\ldots,v_k \in \h^n$ be such that 
\begin{displaymath}
K(u_i,u_j)=K(v_i,v_j)\quad \forall i,j.
\end{displaymath}
Then there exists $g \in \Sp(n)$ such that $g(u_i)=v_i$ for all $i$.
\end{Proposition}

\proof
Let $Q=(q_{ij})=(K(u_i,u_j))$, and denote by $d$ its Moore rank. Then $\mathrm{span}_\h (u_1,\ldots, u_k)$ 
and $\mathrm{span}_\h (v_1,\ldots, v_k)$ have quaternionic dimension $d$.  Without loss of generality, we assume
that $u_1,\ldots,u_d$ are $\h$-linearly independent, or equivalently that
\begin{displaymath}
P=\left(\begin{matrix}q_{11}&\dots&q_{1d}\\\vdots&&\vdots\\q_{d1}&\ldots&q_{dd}\end{matrix}\right)
\end{displaymath}
is invertible. Then $v_1,\ldots,v_d$ are also $\h$-linearly independent. Denoting $P^{-1}=(p^{ij})$,
we have for $r=d+1,\ldots,k$
\begin{align*}
u_r =\sum_{i,j=1}^d u_i p^{ij} q_{jr},\qquad
v_r  =\sum_{i,j=1}^d v_i p^{ij} q_{jr}.
\end{align*}
{If $d=n$, the $\h$-linear map $g$ which sends $u_i$ to $v_i$ preserves $K$ and hence belongs to $\Sp(n)$. If
$d<n$, we may complete $u_1,\ldots,u_d$ (resp. $v_1,\ldots,v_d$) to a basis of $\h^n$ by choosing $K$-orthonormal
vectors in the quaternionic orthogonal complement of $\mathrm{span}_\h (u_1,\ldots,u_d)$  (resp. $\mathrm{span}_\h
(v_1,\ldots,v_d)$). Again, we obtain a map $g \in \Sp(n)$ which maps $u_1,\ldots,u_d$ to
$v_1,\ldots,v_d$.}
\endproof

\begin{Proposition}\label{prop_psi}
Let $V \in \Gr_k$. Denote by $\pi_V:\h^2 \to V$ the orthogonal projection. Given an orthonormal basis
$u_1,\ldots,u_k$ of $V$, we define the endomorphism $\psi_V \in \mathrm{End}(V)$ by                                        
\begin{displaymath}
 \psi_V(y) :=\pi_V \sum_{r=1}^k u_rK(u_r,y)
\end{displaymath}
and set $Q=(q_{ij})_{i,j}:=(K(u_i,u_j))_{i,j}$. Then 
\begin{enumerate}[\rm i)]
 \item $\psi_V$ is independent of the choice of the orthonormal basis $u_1,\ldots,u_k$ of $V$.
\item $\psi_V$ is self-adjoint with respect to the euclidean scalar product on $V$.  
\item If $g \in \spsp$, then $\psi_{gV}=g \circ \psi_V \circ g^{-1}$. In particular, the eigenvalues of $\psi_V$ only
depend on the orbit of $V$. 
\item\label{req2}  The matrix of $\psi_V$ with respect to the basis $u_1,\ldots,u_k$ is $\re Q^2$.
\end{enumerate}
\end{Proposition}
\proof
All claims follow from a straightforward computation.
\endproof

We remark that the endomorphism $\psi_V$ admits the following interpretation:
\begin{displaymath}
 \langle x,\psi_V(y)\rangle =c\int_{\Sp(1)}\langle \pi_V(x\xi),\pi_V(y\xi)\rangle d\xi,\qquad x,y\in V,
\end{displaymath}
where $d\xi$ is the Haar measure on $\Sp(1)$ and $c$ is a non-zero constant.

\subsection{The quotient space $\mathrm{Gr}_2\,/\,\mathrm{Sp}(2)\,\mathrm{Sp}(1)$}

\begin{Theorem} \label{thm_2_orbit}
The quotient $\Gr_2/\spsp$ can be homeomorphically identified with the quotient 
\begin{displaymath}
X_2 := \{\lambda \in [-1,1]\} / \{\pm 1\}
\end{displaymath}
in such a way that $[\lambda]\in X_2$ corresponds to the orbit of
\begin{displaymath}
 V=\spann\{(\cos\theta_1,\sin\theta_1),(\cos\theta_2,\sin\theta_2)\mathbf i\}
\end{displaymath}
with $\lambda=\cos(\theta_1-\theta_2)$.
\end{Theorem}

\proof
Let $V \subset \h^2$ be a two-plane. Choose an orthonormal basis $u_1,u_2$ of $V$. Then $K(u_1,u_2)$ is purely
quaternionic and its norm is bounded by $1$. By using conjugation by an element $\xi \in \Sp(1)$, we may assume that
$K(u_1,u_2)=\lambda \mathbf i$ for some $\lambda \in [-1,1]$. We send the orbit of $V$ to $\lambda$. It is easily checked that
this map is well-defined, a homeomorphism, and fulfills the condition of the statement.
\endproof

\subsection{The quotient space $\mathrm{Gr}_3\,/\,\mathrm{Sp}(2)\,\mathrm{Sp}(1)$}

\begin{Lemma} Under the hypotheses of Proposition \ref{prop_psi} with $k=3$, the following statements are
equivalent:
\begin{enumerate}[\rm i)]
 \item $u_1,u_2,u_3$ is a basis consisting of eigenvectors of $\psi_V$.
\item $q_{12},q_{13},q_{23}$ are pairwise orthogonal in $\im \h$. 
\item $\re Q^2$ is diagonal. 
\end{enumerate}
In this case, the diagonal entries of $\re Q^2$ are the eigenvalues of $\psi_V$. 
\end{Lemma}
\proof This follows easily from claim \ref{req2}) in Proposition \ref{prop_psi}.
\endproof

For each triple $\lambda=(\lambda_{12},\lambda_{13},\lambda_{23}) \in [-1,1]^3$, we denote by $M_\lambda$ the quaternionic $3 \times
3$-matrix 
\begin{displaymath}
 M_\lambda:=\left(\begin{array}{c c c} 1 & \lambda_{12}\i & \lambda_{13}\j \\
             -\lambda_{12}\i & 1 & -\lambda_{23}\k \\
-\lambda_{13}\j & \lambda_{23}\k & 1 
            \end{array}\right).
\end{displaymath}
Let 
\begin{displaymath}
 X_3:=\{\lambda_{pq} \in [-1,1], 1 \leq p < q \leq 3: \rank M_\lambda \leq 2\} /( \mathbb{Z}_2^3 \times
\mathcal{S}_3),
\end{displaymath}
where the action of $\mathbb{Z}_2^3 \times
\mathcal{S}_3$ is given by equations \eqref{actionZ},\eqref{actionS}.

\begin{Proposition}\label{prop_lambda_3}
 Given $V\in \Gr_3$, there is a unique $[\lambda]\in X_3$ such that
 \begin{equation}\label{KqMq}
 K(u_i,u_j)=(M_\lambda)_{i,j},\qquad i,j=1,2,3,
\end{equation}
for some $u_1,u_2,u_3$ spanning an element of the orbit of $V$.
\end{Proposition}

\proof Let $u_1,u_2,u_3\in V$ be an orthonormal basis of eigenvectors of $\psi_V$, and denote $q_{ij}=K(u_i,u_j)$. By the previous lemma, the pure quaternions $q_{12}, q_{13}, q_{23}$ are pairwise orthogonal. Hence there exist  
$\lambda_{12}, \lambda_{13}, {\lambda_{23} \in [-1,1]}$ such that $\lambda_{12} \i,\lambda_{13} \j,$ ${-\lambda_{23}\k
\in \im \h}$ may be mapped to $q_{12},q_{13}, q_{23}$ by a rotation. Let this rotation be $q\mapsto \xi q\bar\xi$ with $\xi\in \Sp(1)$, and let us replace $u_i$ by $u_i\xi$ (without changing the notation). Then, equation \eqref{KqMq} holds. Since $u_1,u_2,u_3$ are
linearly dependent over $\h$, the hyperhermitian matrix $M_\lambda$ has Moore rank at most $2$. Hence
$\lambda=(\lambda_{12},\lambda_{13},\lambda_{23})$ defines a class in  $X_3$. This shows the existence of $[\lambda]$. 

In order to show uniqueness, note that $\re M_\lambda^2$ is diagonal. Hence, by \ref{req2}) of Proposition
\ref{prop_psi}, the orthonormal
basis $u_1,u_2,u_3$ in the statement must consist of eigenvectors of $\psi_V$ (or of $\psi_{gV}$ for some $g\in
\spsp$). 

If $\psi_V$ has three different eigenvalues, then the only freedom
in choosing these vectors is to permute them or to reflect some of them. This results in the action of the group
$\mathbb{Z}_2^3 \times \mathcal{S}_3$ on $\lambda$, so $[\lambda]$ does not depend on the basis. 

If, however, $\psi_V$ has repeated eigenvalues, there are different orthonormal bases consisting of eigenvectors. Let
$u_i, u_i'$ be two such bases, related by $u_i=a_{ij}u'_j$ with $A=(a_{ij})\in \SO(3)$. Take $Q=(K(u_i,u_j))_{i,j}$ and
$Q'=AQA^t=(K(u_i',u_j'))_{i,j}$. We will show that $Q,Q'$ are $\Sp(1)$-conjugate to each other. Hence, the
corresponding
matrices $M_\lambda,M_{\lambda'}$ are $\Sp(1)$-conjugate. It is easy to check that this implies $[\lambda]=[\lambda']$. 

We distinguish two cases depending on the multiplicities of the eigenvalues of
$\psi_V$. 

\textbf{Case 1.} 
Suppose that $\psi_V$ has exactly one double eigenvalue. By reordering the bases, we may assume that the corresponding eigenspace is  $\spann\{u_1,u_2\}=\spann\{u'_1,u'_2\}$, and
\begin{displaymath}
 A=\left(\begin{matrix}\cos \alpha&\sin\alpha&0\\
          -\sin\alpha&\cos\alpha&0\\
	0&0&1
         \end{matrix}\right).
\end{displaymath}
Then $Q'=AQA^t$ has entries $q_{12}'=q_{12}=\lambda_{12}\i$, and
\begin{displaymath}
\left(\begin{matrix}q_{13}'\\
          q_{23}'
\end{matrix}\right)=
\left(
\begin{matrix}\cos \alpha&\sin\alpha\\
          -\sin\alpha&\cos\alpha
\end{matrix}\right)\left(\begin{matrix}\lambda_{13}\j\\
          -\lambda_{23}\k
\end{matrix}\right).
\end{displaymath}
On the other hand, repetition of the eigenvalues means
\begin{displaymath}
 1+\lambda_{12}^2+\lambda_{13}^2=1+\lambda_{12}^2+\lambda_{23}^2
\end{displaymath}
which yields $\lambda_{13}=\epsilon\lambda_{23}$ for some $\epsilon=\pm 1$. Let $\zeta=\cos\frac{\alpha}{2}+\epsilon\sin \frac{\alpha}{2}\i$. Then $Q''=\zeta Q'\bar\zeta$ has entries $q_{12}''=q_{12},q_{13}''=q_{13}, q_{23}''=\epsilon q_{23}$
Since the Moore determinants of $Q,Q''$ vanish, it follows from \eqref{eq_moore3} that $\epsilon=1$ or
$\lambda_{12}\lambda_{13}=0$ or $\lambda_{13},\lambda_{23}$. The latter case can also be reduced to
$\epsilon=1$ by changing the sign of $\lambda_{12},\lambda_{23}$. Hence, $Q',Q$ are
$\Sp(1)$-conjugate to each other, so $[\lambda]=[\lambda']$.

\textbf{Case 2.} Suppose that $\psi_V$ has one triple eigenvalue. Then
\begin{displaymath}
 \lambda_{12}^2+\lambda_{13}^2=\lambda_{12}^2+\lambda_{23}^2=\lambda_{13}^2+\lambda_{23}^2,
\end{displaymath}
so $\lambda_{12}^2=\lambda_{13}^2=\lambda_{23}^2$. By changing signs of $\lambda_{13},\lambda_{23}$, we can assume that $\lambda_{12}=\lambda_{13}$. Then
\begin{displaymath}
 q_{12}'=(a_{11}a_{22}-a_{12}a_{21})\i+(a_{11}a_{23}-a_{13}a_{21})\j+(a_{13}a_{22}-a_{12}a_{23})\k.
\end{displaymath}
Since $A \in \SO(3)$, the wedge product of the first two rows equals the third one, hence 
\begin{displaymath}
 q_{12}'=a_{33}\i-a_{32}\j-a_{31}\k. 
\end{displaymath}
Similarly,
\begin{align*}
 q_{13}'&=-a_{23}\i+a_{22}\j+a_{21}\k,\\
q_{23}'&=a_{13}\i-a_{12}\j-a_{11}\k.\\
\end{align*}
Hence, each $q_{ij}'$ with $i\neq j$ is the image of  $q_{ij}$ under a common rotation of $\R^3\equiv\im\h$.
 Therefore, $Q'$ is an $\Sp(1)$-conjugate of $Q$, and $[\lambda]=[\lambda']$. 

\endproof

\begin{Theorem} \label{thm_3_orbit}
There exists a homeomorphism $X_3 {\cong} \Gr_3/\spsp$ mapping $[\lambda]\in X_3$ to the orbit of a plane
spanned by $v_1, v_2, v_3$ such that
\begin{displaymath}
 K(v_i,v_j)=(M_\lambda)_{i,j},\qquad i,j=1,2,3.
\end{displaymath}
\end{Theorem}

\proof
Given $V \in \Gr_3$, let $[\lambda]\in X_3$ be given by Proposition \ref{prop_lambda_3}. Clearly $[\lambda]$
only depends on the $\spsp$-orbit of $V$ in $\Gr_3$. Hence, $V \mapsto [\lambda]$
defines a map {$\Phi:\Gr_3/\spsp \rightarrow X_3$}. 

Let us show that $\Phi$ is bijective. To show injectivity, suppose that $U,V \in  \Gr_3$ are mapped to
the same $[\lambda]\in X_3$. This means that $U$ and $V$ admit respective bases  $u_1,u_2,u_3$ and $v_1,v_2,v_3$, such
that
\begin{displaymath}
K(u_i\zeta,u_j\zeta)=K(v_i\xi,v_j\xi)=M_\lambda
\end{displaymath}
for certain $\zeta,\xi\in\Sp(1)$. By
Proposition \ref{prop_inj}, there exists $g\in\Sp(2)$ such that $g(u_i \zeta)=v_i \xi$. Hence $V=g(U)\zeta \bar \xi$,
so 
$U$ and $V$ belong to the same $\spsp$-orbit.

To see surjectivity, it is enough to apply Proposition \ref{prop_surj} with $Q=M_\lambda$. 

Since $\Gr_3$ is compact and $X_3$ is Hausdorff, it remains only to prove that $\Phi$ is continuous. 

{Let $(V^m)$ be a sequence of $3$-planes converging to the $3$-plane $V$ in $\Gr_3$. Let $(u_1^m, u_2^m,
u_3^m)$
be an orthonormal basis of $V^m$ and $\lambda^m=(\lambda_{12}^m,\lambda_{13}^m,\lambda_{23}^m)$ as in
Proposition \ref{prop_lambda_3}. By compactness, there exists a subsequence $m_1,m_2,\ldots$ such that
$(u_1^{m_l},u_2^{m_l},u_3^{m_l})$ converges to an orthonormal basis $(u_1,u_2,u_3)$ of $V$. Hence $\lambda^{m_l} \to
\lambda$ for some $\lambda=(\lambda_{12},\lambda_{13},\lambda_{23})$. Then $\Phi(V)=[\lambda]$ and it follows that
$\Phi(V_{m^l})$ converges to $\Phi(V)$. 

Since we may apply the same argument to any subsequence of a given sequence, we obtain the following: every subsequence of $(V_m)$ contains a subsequence such that the images under $\Phi$ converge to $\Phi(V)$. But this implies that the images under $\Phi$ of the original sequence converge to $\Phi(V)$.}      
\endproof

\begin{Corollary}\label{coro_angles_3}
Given $[\lambda]\in X_3$, there exist $\theta_1,\theta_2,\theta_3$ such that $$\lambda_{ij}=\cos(\theta_i-\theta_j),$$ and the orbit corresponding to $[\lambda]$ contains the plane
 \begin{displaymath}
 V=\spann\{(\cos\theta_1,\sin\theta_1),(\cos\theta_2 ,\sin\theta_2)\mathbf i, (\cos\theta_3 ,\sin\theta_3)\mathbf j\}.
\end{displaymath}
\end{Corollary}

\proof
By Theorem \ref{thm_3_orbit}, the orbit corresponding to $[\lambda]$ contains a plane $V$ admitting an
orthonormal basis $v_1,v_2,v_3$ such that $K(v_i,v_j)=(M_\lambda)_{i,j}$. Since $\Sp(2)$ acts transitively on the unit
sphere of $\h^2$, we can assume $v_1=(1,0)$. From $K(v_1,v_2)=\lambda_{12}\mathbf i$, we deduce that
$v_2=(\lambda_{12}\mathbf i,w)$ for some $w\in \h$. By applying an element of  $\Sp(1)$ to the second component of
$\h^2$, we may assume that $w$ and $\mathbf i$ are parallel, $w\|\mathbf i$.
Together with $K(v_2,v_3)=\lambda_{23}\mathbf j$, this implies that $v_3=(a\mathbf j, b\mathbf j)$ for some $a,b\in\R$. Therefore, $V$ agrees with the given description.
\endproof

\subsection{The quotient space $\mathrm{Gr}_4\,/\,\mathrm{Sp}(2)\,\mathrm{Sp}(1)$}\mbox{}

Let $V \subset \h^2$ be a $4$-plane. Given an orthonormal basis $u_1,\ldots,u_4$ of $V$, we set
$Q:=(K(u_p,u_q))_{p,q}$.
Clearly the Moore rank of $Q$ is at most $2$ and $\tr Q=4$. We call $V$ {\it degenerated} if $Q$ has Moore
eigenvalues $(2,2,0,0)$ and {\it non-degenerated} otherwise. This notion is independent of the choice of the
orthonormal basis. 

Note that if $\re Q^2=2 \Id$ (which is equivalent to $\psi_V=2 \Id$), then $Q$ is degenerated. Indeed,
if $\lambda,4-\lambda$ are the non-zero Moore
eigenvalues of $Q$, then $\lambda^2+(4-\lambda)^2=\tr Q^2=8$ which implies that $\lambda=2$.  

\begin{Lemma} \label{lemma_density}
Non-degenerated planes are dense in $\Gr_4$.   
\end{Lemma}

\proof
Consider the continuous map which sends $g \in \SO(8)$ to the plane $V$ spanned by the first four columns in $\R^8 \cong \h^2$.
Let $u_1,\ldots,u_8$ be the columns of $g$ and $Q:=(K(u_p,u_q))_{p,q}$. Then $V$ is non-degenerated if
and
only if $\tr Q^2 \neq 8$. Clearly the function $\tr Q^2-8$ is a polynomial
function on the irreducible algebraic variety $\SO(8)$. Since this function does not vanish identically
on $\SO(8)$, its zero set does not contain any
open set. 
\endproof

\begin{Proposition} \label{prop_claim}
In each $\spsp$-orbit of $\Gr_4$ there is an element with an orthonormal basis $v_1, v_2, v_3,v_4$ such
that
each $v_i=(v_{i1},v_{i2})\in\h^2$ has parallel components; i.e. $v_{i1}\|v_{i2}$ as vectors of $\h\equiv\R^4$ for $i=1,\ldots, 4$.
\end{Proposition}

\proof
By Lemma \ref{lemma_density}, non-degenerated $4$-planes are dense in $\Gr_4$. By continuity it
is enough to prove the statement for non-degenerated planes.  

Let $V \in \Gr_4$ be non-degenerated and let $u_1,\ldots,u_4$ be a basis consisting of eigenvectors of
$\psi_V$. 
Define 
\begin{displaymath}
Q:=(K(u_m,u_l))_{m,l=1,\ldots,4}. 
\end{displaymath}

Since $u_1,\ldots,u_4$ are eigenvectors of $\psi_V$, the matrix $\re Q^2$ is diagonal.
Moreover, $\tr Q=4$ and the Moore rank of $Q$ is at most $2$. We can therefore write $Q=A^*DA$, where $A=(a_{ij}) \in
\Sp(2)$ and
$D= \diag(\delta,4-\delta,0,0), \delta \in [0,4]$. Since $V$ is non-degenerated, we have $ \delta \neq 2$, hence $\re
Q^2 \neq 2 \Id$. 

We claim that $a_{1m}, m=1,\ldots,4$ are pairwise orthogonal in $\mathbb{H}$, and the same holds for
$a_{2m},m=1,\ldots,4$. 
For instance, we have 
\begin{displaymath}
 q_{12}=\delta \bar a_{11}a_{12}+(4-\delta) \bar a_{21} a_{22}
\end{displaymath}
and 
\begin{displaymath}
 (Q^2)_{12}= \delta^2 \bar a_{11}a_{12}+(4-\delta)^2 \bar a_{21} a_{22}.
\end{displaymath}
The real part of these two quaternions vanishes if and only if $\bar a_{11}a_{12}$ and $\bar a_{21} a_{22}$ are pure
quaternions (here we use that $\delta \neq 2$). 

The matrix $A$ can be left multiplied by a diagonal matrix with entries in $\Sp(1)$ and $Q$ remains unchanged.
Since this
action is transitive on the unit sphere in each summand of $\h^2=\h \oplus \h$, we can
assume that $a_{14},a_{24}\in\R^+$. Also, we can conjugate $A$ by an element $\xi \in \Sp(1)$. 
The effect is that also $Q$ is conjugated by $\xi$, which is equivalent to multiplying $V$ by $\xi$ from the right.

The vectors $ (\sqrt{\delta}a_{1m},\sqrt{4-\delta}a_{2m}), m=1,\ldots,4$ form an orthonormal basis of a
$4$-plane in the same orbit as $V$. We may therefore assume that $V$ is spanned by the vectors 
\begin{align}\label{u1}
 u_1 & =(\sqrt{\delta}a_{11}, \sqrt{4-\delta}a_{21}) =:(\cos\theta_1\,\i,\sin\theta_1\,w_1),\\
u_2 & =(\sqrt{\delta}a_{12}, \sqrt{4-\delta}a_{22}) =:(\cos\theta_2\,\j,\sin\theta_2\,w_2),\\
u_3 & =(\sqrt{\delta}a_{13}, \sqrt{4-\delta}a_{23}) =:(\cos\theta_3\,\k,\sin\theta_3\,w_3),\\\label{u4}
 u_4 & =(\sqrt{\delta}a_{14}, \sqrt{4-\delta}a_{24}) =:(\cos\theta_4,\sin\theta_4),
\end{align}
where $w_1,w_2,w_3$ is an orthonormal basis of $\R^3\equiv \im \mathbb H$.  

By changing the sign of some $u_m, w_m$
we can suppose that $0 \leq \theta_{1},\ldots,\theta_{4}\leq \frac\pi 2$.

Since $A \in \Sp(2)$, we have $\sum \bar a_{1m} a_{2m}=0$, i.e. 
\begin{equation}\label{condition}
 \sin (2\theta_4)-\sin(2\theta_1)\, \i\cdot w_1 -\sin(2\theta_2)\, \j\cdot
w_2-\sin (2\theta_3)\, \k\cdot w_3=0.
\end{equation}

Considering the imaginary part we deduce
\begin{displaymath}
 \sin (2\theta_m) w_{mn}=\sin (2\theta_n) w_{nm},\qquad m,n=1,2,3,
\end{displaymath}
where $w_{mn}$ are the coordinates of $w_m$ with respect to the basis $\i,\j,\k$ of $\R^3$; i.e, the matrix
$M=(\sin (2\theta_m) w_{mn})_{m,n=1,2,3}$ is symmetric. Let $d_m:=\sin 2\theta_m$, 
$D:=\diag(d_1,d_2,d_3)$ and
$O:=(w_1,w_2,w_3) \in \mathrm{O}(3)$. Then $M=DO$ and hence $DO=O^tD, OD=DO^t$. Therefore $OD^2=DO^tD=D^2O$, i.e.
\begin{displaymath}
 (d_i^2-d_j^2)o_{ij}=0.
\end{displaymath}

We consider three cases according to the multiplicities of the entries in $D$.

\textbf{Case 1.}
If $\#\{d_i\}=3$ then $O$ is diagonal and the statement is trivial.

\textbf{Case 2.} $\#\{d_i\}=2$ and $O$ contains a row with zeros outside the diagonal
position, i.e. up to a simultaneous reordering of rows and columns, $D$ and $O$ have the form
\begin{displaymath}
 D=\left(\begin{matrix} d_1 & 0 & 0\\0 & d_1 & 0\\ 0 & 0 &d_3\end{matrix}\right),\qquad O=\left(\begin{matrix}
\cos\alpha&\sin\alpha&0\\\sin\alpha&-\cos\alpha&0\\0&0&\varepsilon\end{matrix}\right), \quad\varepsilon=\pm 1.
\end{displaymath}

After reordering $u_1,u_2,u_{3}$ and conjugating by a suitable element of $\Sp(1)$ we have
\begin{align*}
u_1 &=(\cos\theta_1\, \i,\sin\theta_1(\cos\alpha\, \i+\sin\alpha\, \j))\\
u_2 &=(\cos\theta_2\, \j,\sin\theta_2(\sin\alpha\, \i-\cos\alpha\, \j))\\
u_3 &=(\cos\theta_3\, \k,\varepsilon\sin\theta_3\, \k)\\
u_4 &=(\cos\theta_4,\sin\theta_4)
\end{align*}
with $\sin 2\theta_{1}=\sin 2\theta_{2}$. Thus, either $\theta_{2}=\theta_{1}$ or
$\theta_{2}=\frac{\pi}{2}-\theta_{1}$. 

By considering the real part of \eqref{condition} we deduce $\sin 2\theta_3=\sin 2\theta_4$ and
$\varepsilon=-1$. 

We consider three cases. 
\begin{itemize}
\item If $\theta_{2}=\theta_{1}$, we set $u_1':=\cos \frac{\alpha}{2} u_1+\sin
\frac{\alpha}{2}
u_2, u_2':=-\sin \frac{\alpha}{2} u_1+\cos \frac{\alpha}{2}
u_2, u_3'=u_3, u_4'=u_4$. Then, the first and second components of $u_i'\in  \h^2$ are
parallel for each $1\leq i\leq 4$. 

\item If $\theta_{3}=\theta_{4}$, we set $u_1':=u_1, u_2':=u_2, u_3':=\cos \frac{\alpha}{2} u_3+\sin \frac{\alpha}{2}
u_3, u_4':=-\sin \frac{\alpha}{2} u_3+\cos \frac{\alpha}{2} u_4$. Again we obtain an orthonormal basis of $V$ that
satisfies the statement.
 \item If $\theta_{2}=\frac{\pi}{2}-\theta_{1}$ and
$\theta_{4}=\frac{\pi}{2}-\theta_{3}$, then one checks that $\re(Q^2)=2\, \Id$, contradicting our assumption.
\end{itemize}

\textbf{Case 3.} $D$ is a multiple of the identity. 

Then $\sin 2\theta_{m}=c \neq 0$ for $m=1,2,3$. The real part of
\eqref{condition} is
\begin{displaymath}
\sin 2\theta_4+c\,\mathrm{tr}O=0.
\end{displaymath}
Since $O$ is orthogonal and diagonalizable, it
has eigenvalues $1,1,1$ or $1,1,-1$ or $1,-1,-1$ or $-1,-1,-1$. In the first and last cases, $O$ is diagonal and we are
done. Otherwise $\tr O=\pm
1$. Since $\sin 2\theta_{m}\geq 0$, we deduce that $\tr  O=- 1$, i.e. $O$ has eigenvalues $1,-1,-1$, and $\sin
2\theta_{4}=c$.  

Therefore every two angles  $\theta_{m},\theta_{n}, {1\leq m,n \leq 4}$ are equal or complementary. If
$\theta_1,\ldots,\theta_{4}$ contain
exactly two pairs of equal angles, then one checks that $\re(Q^2)=2\, \Id$, again contradicting our assumption.
Hence at least three angles $\theta_{m}$ are equal. By reordering, we may assume that 
$\theta_{1}=\theta_{2}=\theta_{3}$. Then we write
\begin{displaymath}
O=P^t\left(\begin{matrix} -1&0&0\\ 0&-1&0\\0&0&1\end{matrix}\right)P,
\end{displaymath}
where $P \in \mathrm{O}(3)$ and set 
\begin{displaymath}
\left(\begin{matrix} u_1'\\ u_2'\\ u_3'\end{matrix}\right):=P\left(\begin{matrix}u_{1} \\ u_{2}\\
u_{3}\end{matrix}\right) \quad u_4':=u_4.
\end{displaymath}
Then, the first and second components of each $u_i'$ are parallel vectors in $\h$. 
\endproof

\begin{Corollary}
Every $V\in \Gr_4$ admits an orthonormal basis $u_1,\ldots, u_4$ such that $q_{ij}=K(u_i,u_j)$  satisfy
\begin{itemize}
 \item $q_{12},q_{13},q_{23}$ are pairwise orthonormal
\item $q_{12}\|q_{34}, q_{13}\|q_{24}, q_{14}\|q_{23}$.
\end{itemize}\end{Corollary}
\proof
It is enough to check the statement for one plane in each $\spsp$-orbit of $\Gr_4$. By the previous proposition, we may
assume that $V$ admits an orthonormal basis $u_1,\ldots, u_4$ with $u_{i1},u_{i2}$ both parallel to some
$\xi_i\in\h\setminus\{0\}$ for each $i$. 
Since $u_1,\ldots, u_4$ are orthogonal, so are $\xi_1,\ldots,\xi_4$. Since $q_{ij}\|\bar\xi_i\xi_j$, we get $q_{ij}\bot
q_{ik}$ if $j\neq k$. The statement follows.
\endproof

Given $\lambda_{pq} \in [-1,1], 1 \leq p < q \leq 4$, we define the quaternionic matrix 
\begin{displaymath}
 M_\lambda:=\left(\begin{array}{c c c c} 1 & \lambda_{12}\i & \lambda_{13}\j & \lambda_{14}\k \\
             -\lambda_{12}\i & 1 & -\lambda_{23}\k & \lambda_{24}\j\\
-\lambda_{13}\j & \lambda_{23}\k & 1 & -\lambda_{34}\i\\
-\lambda_{14}\k & -\lambda_{24}\j & \lambda_{34}\i & 1
            \end{array}\right).
\end{displaymath}
Let 
\begin{displaymath}
 X_4:=\{\lambda_{pq} \in [-1,1], 1 \leq p < q \leq 4: \rank M_\lambda \leq 2\} / (\mathbb{Z}_2^4 \times
\mathcal{S}_4),
\end{displaymath}
where the action of $\mathbb{Z}_2^4 \times
\mathcal{S}_4$ is given by equations \eqref{actionZ},\eqref{actionS}. 

\begin{Proposition}
 Given $V\in \Gr_4$, there is a unique $[\lambda]\in X_4$ such that
 \begin{displaymath}
 K(u_i,u_j)= (M_\lambda)_{i,j},\qquad i,j=1,2,3,4,
\end{displaymath}
for some $u_1,\ldots,u_4$ spanning an element of the orbit of $V$.
\end{Proposition}

\proof Let $u_1,\ldots, u_4$ be given by the previous corollary. 
Using a rotation $q \mapsto \xi q \bar \xi$, we may map $q_{12}$ to a multiple of $\i$, $q_{13}$ to a multiple of $\j$
and $q_{14}$ to a multiple of $\k$. For $i=1,\ldots, 4$ take $u_i\xi$ and denote it again by $u_i$. Then,  
\begin{align}
 K(u_1,u_2) & = \lambda_{12} \i \label{eq_k12}\\
K(u_1,u_3) & = \lambda_{13} \j \label{eq_k13}\\
K(u_1,u_4) & = \lambda_{14} \k \label{eq_k14}\\
K(u_2,u_3) & = -\lambda_{23} \k \label{eq_k23}\\
K(u_2,u_4) & = \lambda_{24} \j \label{eq_k24}\\
K(u_3,u_4) & = -\lambda_{34} \i \label{eq_k34}
\end{align}
for real numbers $\lambda_{pq} \in [-1,1], 1 \leq p<q \leq 4$.
Since any $3$ vectors in $\h^2$ are linearly dependent over $\h$, the rank of the matrix $Q:=M_\lambda$ is at most $2$. This shows the existence part of the statement.

In order to prove uniqueness, let $A=(a_{ij})\in \SO(4)$ and suppose that  $u_i'=a_{ij}u_j$ is
another basis of
$V$ such that $Q'=AQA^t$ is $\Sp(1)$-conjugate to $(M_{\lambda'})_{ij}$ for some $[\lambda']\in X_4$. Then $\re Q^2,\re
(Q')^2$ are both diagonal. By Proposition \ref{prop_psi}, the orthonormal
bases $u_1,\ldots,u_4$ and $u'_1,\ldots, u'_4$ consist both of eigenvectors of $\psi_V$.  We need to show that $Q,Q'$
are $\Sp(1)$-conjugates of each other, {which will imply that} $[\lambda]=[\lambda']$.

If $\psi_V$ has no multiple eigenvalues, then the two bases coincide up to signs and order. Hence $[\lambda]=[\lambda']$.

Next we consider different cases according to the multiplicities of the eigenvalues of
$\psi_V$. 

\textbf{Case 1.} 
Suppose that $\psi_V$ has exactly one double eigenvalue. By reordering the bases, we may assume that the corresponding
eigenspace is  $\spann\{u_1,u_2\}=\spann\{u'_1,u'_2\}$, and
\begin{displaymath}
 A=\left(\begin{matrix}\cos \alpha&\sin\alpha&0&0\\
          -\sin\alpha&\cos\alpha&0&0\\
	0&0&1&0\\
0&0&0&1
         \end{matrix}\right).
\end{displaymath}
Then $Q'=AQA^t$ has entries $q_{12}'={q_{12}}$, $q_{34}'={q_{34}}$, and
\begin{displaymath}
\left(\begin{matrix}q_{13}'&q_{14}'\\
          q_{23}'&q_{24}'
\end{matrix}\right)=
\left(
\begin{matrix}\cos \alpha&\sin\alpha\\
          -\sin\alpha&\cos\alpha
\end{matrix}\right)\left(\begin{matrix}\lambda_{13}\j&\lambda_{14}\k\\
          -\lambda_{23}\k&\lambda_{24}\j
\end{matrix}\right).
\end{displaymath}

Our assumption is that each row  and each column in $Q'$ has orthogonal entries. This implies that either
$\sin\alpha\cos\alpha=0$, in which case everything follows trivially, or $\lambda_{13}=\epsilon\lambda_{23},
\lambda_{14}=\epsilon\lambda_{24}$ for some $\epsilon=\pm 1$. Since the $3\times 3$ upper left minors of $Q,Q'$
vanish,
we have $\epsilon=1$ (except if $\lambda_{13}\lambda_{23}=0$, in which case we may choose $\epsilon=1$ as well). It
follows that $Q'=\bar \zeta Q\zeta$ with $\zeta=\cos\frac{\alpha}{2}+\sin
\frac{\alpha}{2}\i$. 

\textbf{Case 2.} 
Suppose that $\psi_V$ has two different double eigenvalues. We may assume that $A$ has the form
\begin{displaymath}
 A=\left(\begin{matrix}\cos \alpha&\sin\alpha&0&0\\
          -\sin\alpha&\cos\alpha&0&0\\
	0&0&\cos\beta&\sin\beta\\
0&0&-\sin\beta&\cos\beta
         \end{matrix}\right).
\end{displaymath}
Then $\lambda_{13}^2+\lambda_{14}^2=\lambda_{23}^2+\lambda_{24}^2$ as well as
$\lambda_{13}^2+\lambda_{23}^2=\lambda_{14}^2+\lambda_{24}^2$, which implies that $\lambda_{13}^2=\lambda_{24}^2$ and
$\lambda_{14}^2=\lambda_{23}^2$. 

By changing some  sign  if necessary, we may assume that $\lambda_{13}=\lambda_{24}$. The rank $2$ condition of $Q$
leads to $\lambda_{14}=\lambda_{23}$ or $\lambda_{13}\lambda_{14}=0$ or $\lambda_{12}=\lambda_{34}=0$. The third
possibility is excluded by the assumption that the eigenvalues are different, and the second one also allows to suppose
$\lambda_{14}=\lambda_{23}$.

The upper right
square of $Q$  is thus given by 
\begin{displaymath}\left(\begin{array}{c c} q_{13} & q_{14} \\ q_{23}  &
q_{24} \end{array}\right)=
 \lambda \left(\begin{array}{c c} \cos(\theta) \j & \sin(\theta) \k \\ -\sin(\theta)\k  &
\cos(\theta) \j \end{array}\right),
\end{displaymath}
where $\lambda:=\sqrt{\lambda_{13}^2+\lambda_{14}^2}$. The upper right square of $Q'$ is 
{\footnotesize \begin{displaymath} 
 \lambda \left(\begin{array}{c c} \cos(\theta)\cos(\alpha-\beta) \j-\sin(\theta) \sin(\alpha-\beta)\k &
\cos(\theta)\sin(\alpha-\beta) \j+\sin(\theta)\cos(\alpha-\beta)\k \\ -\cos(\theta)\sin(\alpha-\beta)
\j-\sin(\theta)\cos(\alpha-\beta)\k &
\cos(\theta) \cos(\alpha-\beta) \j -\sin(\theta) \sin(\alpha-\beta) \k \end{array}\right).
\end{displaymath}}
The assumption that rows and columns have orthogonal entries implies that either $2\alpha-2\beta$ is a
multiple of $\pi$, or $\sin^2 \theta=\cos^2 \theta$. In the first case, one checks easily that $Q'$ is related to $Q$ by
an element of $\mathbb Z_2^4\times\mathcal S_4$.

Next, suppose that $\sin^2 \theta=\cos^2 \theta=\frac12$. 
In this case $Q$ and $Q'$
differ only by a rotation in the plane $\spann\{\j,\k\}$.  

\textbf{Case 3.} 
Suppose that $\psi_V$ has a triple eigenvalue, say corresponding to the first three vectors of each
basis. Then $A\in  \SO(3)\subset \SO(4)$, and
\begin{displaymath}
 \lambda_{12}^2+\lambda_{13}^2+\lambda_{14}^2= \lambda_{12}^2+\lambda_{23}^2+\lambda_{24}^2=
\lambda_{13}^2+\lambda_{23}^2+\lambda_{34}^2.
\end{displaymath}

Putting $P=(q_{14},q_{24},q_{34})^t=(\lambda_{14}\k,\lambda_{24}\j,-\lambda_{34}\i)^t$ we have
\begin{displaymath}
 PP^*=\left(\begin{matrix} \lambda_{14}^2 & 0 & 0\\
             0 & \lambda_{24}^2 & 0\\
             0 & 0 & \lambda_{34}^2
            \end{matrix}\right)=:D.
\end{displaymath}
By assumption,  $P'=(q_{14}',q_{24}',q_{34}')^t$ has orthogonal entries. Since  $P'=AP$ we deduce that
$D':=P'(P')^*=ADA^t$ is diagonal. After multiplication of $A$ by a permutation matrix, we can assume $D'=D$.

From $AD=DA$ we get three possibilities: either $\lambda_{14}^2,\lambda_{24}^2,\lambda_{34}^2$ has no repetitions and $A$ is the identity, or $\#\{\lambda_{14}^2,\lambda_{24}^2,\lambda_{34}^2\}=2$ and $A$ is a rotation in some
$2$-plane (this case can be handled as Case 1), or $\lambda_{14},\lambda_{24},\lambda_{34}$ have the same
absolute value $\mu$. From the equations above it follows that $\lambda_{12},\lambda_{13},\lambda_{23}$ also have the
same absolute value $\tau$. We may assume that $\lambda_{12},\lambda_{13},\lambda_{14} \geq 0$. Then $\lambda_{23}=\pm
\tau, \lambda_{24}=\pm \mu,\lambda_{34}=\pm \mu$.

Since the upper $3 \times 3$ minor of $Q$ must vanish, {we obtain from \eqref{eq_moore3}} that $\tau \in
\left\{\pm 1, \pm \frac12\right\}$. Checking all possible combinations,
the only matrices of this type of rank $2$ are
 \begin{displaymath}
 Q=\left(\begin{array}{c c c c} 1 & \i & \j & \mu \k \\
             -\i & 1 & -\k & \mu \j\\
-\j & \k & 1 & -\mu \i\\
-\mu \k & -\mu \j & \mu \i & 1
            \end{array}\right),
\end{displaymath}
where $\mu$ is arbitrary. The rest of the proof in this case is analogous to Case 2 in the proof of Proposition \ref{prop_lambda_3}.
 
\textbf{Case 4.} 
Suppose that all eigenvalues of $\psi_V$ are the same. Then 
 \begin{displaymath}
\lambda_{12}^2+\lambda_{13}^2+\lambda_{14}^2= \lambda_{12}^2+\lambda_{23}^2+\lambda_{24}^2=
\lambda_{13}^2+\lambda_{23}^2+\lambda_{34}^2=\lambda_{14}^2+\lambda_{24}^2+\lambda_{34}^2,
\end{displaymath}
which implies that $\lambda_{23}=\epsilon_1 \lambda_{14},
\lambda_{24}=\epsilon_2 \lambda_{13}, \lambda_{34} = \epsilon_3 \lambda_{12}$ with
$\epsilon=(\epsilon_1,\epsilon_2,\epsilon_3) \in \{\pm 1\}^3$. Using the fact that $Q$ has Moore rank $2$ yields two
possibilities
\begin{itemize}
 \item[i)] $\epsilon_1=\epsilon_2=\epsilon_3$
\item[ii)] $\lambda_{12}\lambda_{13}\lambda_{14}=0$. 
\end{itemize}

In case $\rm i)$, we can assume
\begin{displaymath}
Q=\left(\begin{array}{c c c c} 1 & q_{12} & q_{13} & q_{14} \\ -q_{12} & 1 & - q_{14} &  q_{13}\\ -q_{13} &
 q_{14} & 1 & - q_{12}\\-q_{14} &  -q_{13} &  q_{12} & 1
  \end{array}\right).
\end{displaymath}
The conjugation of a matrix of this form by $A\in \SO(4)$ can be described as follows. Let $\Lambda_-^2\mathbb
R^4$ be
the $(-1)$-eigenspace of the Hodge operator $*:\Lambda^2 \R^4 \to \Lambda^2\R^4$. We identify $\Lambda_-^2 \R^4$ with
$\R^3$ by choosing the orthonormal basis ${e_1 \wedge e_2}
- e_3 \wedge e_4, e_1 \wedge e_3 + e_2 \wedge e_4, e_1 \wedge e_4 - e_2 \wedge e_3$. The action of $\SO(4)$ on
$\Lambda^2\mathbb R^4$ preserves $\Lambda_-^2 \R^4\cong \R^3$, which yields a map $\rho:\SO(4)\to \SO(3)$. 

Now consider real $4\times 4$-matrices of the form 
\begin{displaymath}
P:=\left(\begin{array}{c c c c} 1 & x_{12} & x_{13} & x_{14} \\ -x_{12} & 1 & x_{23} & x_{24}\\ -x_{13} & -x_{23} & 1 &
x_{34}\\-x_{14} & -x_{24} & -x_{34} & 1
  \end{array}\right)
\end{displaymath}
and set $\iota(P):=\sum_{1\leq i < j \leq 4} x_{ij} e_i \wedge e_j \in \Lambda^2 \R^4$. Then
$\iota(P) \in \Lambda^2_- \R^4$ if and only if $x_{34}=-x_{12}, x_{24}=x_{13}, x_{23}=-x_{14}$. In this case, 
$\iota(APA^t)=\rho(A)(\iota(P))$ for $A\in \SO(4)$.

Tensorizing everything with $\R^3=\im\h$ we conclude that $Q'=AQA^t$ has the same form as $Q$ and 
\begin{displaymath}
 \left(\begin{matrix}q_{12}'\\q_{13}'\\q_{14}'\end{matrix}\right)=\rho(A)\left(\begin{matrix}q_{12}\\q_{13}\\q_{14}\end{matrix}\right).
\end{displaymath}
Hence, $Q'$ is obtained by applying a rotation of $\mathbb R^3$ to the purely quaternionic coefficients of $Q$; i.e. $Q$ and $Q'$ are $\Sp(1)$-conjugates of each other. 

In case $\rm ii)$, after reordering indices we may suppose $\lambda_{12}=\lambda_{34}=0$. From the rank $2$ condition we
also have
\begin{displaymath}
 \lambda_{13}^2+\lambda_{14}^2=1,\qquad (\epsilon_1\lambda_{14}^2-\epsilon_2\lambda_{13}^2)^2=1.
\end{displaymath}
Hence, $\lambda_{13}=\cos\theta, \lambda_{14}=\sin\theta$ for some $\theta$. Moreover, the second equation yields
$\epsilon_1\epsilon_2=-1$ or $\sin\theta\cos\theta=0$. In both cases, after the action of $\mathbb Z_2^4$ we can assume
$\lambda_{13}=\lambda_{24}=\cos\theta$ and $\lambda_{14}=-\lambda_{23}=\sin\theta$. The matrix $M_\lambda$ is then given
by 
\begin{displaymath}
 M_\lambda=\left(\begin{array}{c c c c} 1 & 0 & \cos \theta \j & \sin \theta \k \\
            0 & 1 & \sin \theta \k & \cos \theta\j\\
-\cos \theta \j & -\sin \theta \k & 1 & 0\\
-\sin \theta \k & -\cos \theta \j & 0 & 1
            \end{array}\right).
\end{displaymath}
Up to permutations, $M_{\lambda'}$ has the same form possibly with a different $\theta$.

The function 
\begin{displaymath}
W \mapsto \min_{u \in W , \|u\|=1}\max_{\xi \in S^3 \cap \im \h} |\pi_W( u\cdot\xi)|
\end{displaymath}
is a $\spsp$-invariant function on $\Gr_4$. It is easily checked that it assumes the value
$\max\{|\cos \theta|,|\sin \theta|\}$ on the plane $V$. The proof is completed by noting that the equivalence class
of $[\lambda]$ only depends on $\max\{|\cos \theta|,|\sin \theta|\}$. 
\endproof

\begin{Theorem} \label{thm_4_orbit}
There exists a homeomorphism $X_4 {\cong}  \Gr_4/\spsp$ mapping $[\lambda]\in X_4$ to the
orbit of a plane spanned by $v_1,\ldots, v_k$ such that
\begin{displaymath}
 K(v_i,v_j)=(M_\lambda)_{i,j},\qquad i,j=1,\ldots,4.
\end{displaymath}
\end{Theorem}
The proof is exactly as in Theorem \ref{thm_3_orbit}.

\begin{Corollary}\label{coro_angles_4}
Given $[\lambda]\in X_k$, there exist $\theta_1,\ldots,\theta_4$ such that $$\lambda_{ij}=\cos(\theta_i-\theta_j),$$ and
the orbit corresponding to $[\lambda]$ contains the plane
 \begin{displaymath}
 V=\spann\{(\cos\theta_1,\sin\theta_1),(\cos\theta_2,\sin\theta_2) \mathbf i, (\cos\theta_3 ,\sin\theta_3)\mathbf
j,(\cos\theta_4,\sin\theta_4)\mathbf k\}.
\end{displaymath}
\end{Corollary}
The proof is analogous to that of  Corollary \ref{coro_angles_3}.

\section{Irreducible representations of $\mathrm{SO}(n)$}

It is well-known that equivalence classes of complex irreducible (finite-dimensional) representations of $\SO(n)$ are
indexed by their highest weights. The possible highest weights are tuples
$\left(\lambda_1,\lambda_2,\ldots,\lambda_{\left\lfloor
\frac{n}{2}\right\rfloor}\right)$ of integers such that  
\begin{enumerate}[\rm i)]
 \item
$\lambda_1 \geq \lambda_2 \geq \ldots \geq \lambda_{\lfloor
\frac{n}{2}\rfloor} \geq 0$ if  $ n $ is odd,
\item $\lambda_1 \geq \lambda_2 \geq \ldots \geq |\lambda_{\frac{n}{2}}| \geq 0$ if $ n $ is  even.
\end{enumerate}

We will write $\Gamma_\lambda$ for any isomorphic copy of an irreducible representation with highest weight $\lambda$. As in \cite{alesker_bernig_schuster}, if $n$ is even and $\lambda=(\lambda_1,\lambda_2,\ldots,\lambda_\frac{n}{2})$ then we set $\lambda':=(\lambda_1,\lambda_2,\ldots,-\lambda_\frac{n}{2})$. It will be useful to use the following notation:  
\begin{displaymath}
 \tilde \Gamma_\lambda:=
\begin{cases} 
\Gamma_\lambda & n \text{ odd or } \lambda_{\frac{n}{2}}=0\\
\Gamma_\lambda \oplus \Gamma_{\lambda'} & n \text{ even and } \lambda_\frac{n}{2} \neq 0.
\end{cases}
\end{displaymath}

The following proposition is well-known, compare \cite{strichartz75, takeuchi73} and  (\cite{olafsson_pasquale}, Lemma 5.3). 

\begin{Proposition}
 Let $\Gr_k(\R^n)$ denote the Grassmann manifold consisting of all $k$-dimensional subspaces in $\R^n$. The
$\SO(n)$-module $L^2(\Gr_k(\R^n))$ decomposes as 
\begin{displaymath}
 L^2(\Gr_k(\R^n)) \cong \bigoplus_\lambda \Gamma_\lambda,
\end{displaymath}
where $\lambda$ ranges over all highest weights such that $\lambda_i=0$ for $i
> \min\{k,n-k\}$ and such that all $\lambda_i$ are even. In particular, it is multiplicity-free.
\end{Proposition}

Let $\Gamma_\lambda$ be an irreducible representation of $\SO(n)$ appearing in $L^2(\Gr_k(\R^n))$. By Schur's
lemma, the
Laplacian $\Delta$ acts by multiplication by some scalar, which was computed by James-Constantine
\cite{constantine_james74}. We will follow the convention $\Delta f:=-\mathrm{div}\circ\nabla f$.

\begin{Proposition} \label{prop_eigenvalues}
 The Laplace-Beltrami operator $\Delta$ of $\Gr_k(\R^n)$ acts on $\Gamma_\lambda$ by the scalar 
\begin{displaymath}
 \sum_{i=1}^{\left\lfloor \frac{n}{2} \right\rfloor} \lambda_i(\lambda_i-2i+n).
\end{displaymath}
\end{Proposition}

We will also need the decomposition of $\Val_k$ as a sum of irreducible $\SO(n)$-modules, which was
obtained recently in \cite{alesker_bernig_schuster}. 
\begin{Proposition}\label{prop_val_decomposition}
The $\SO(n)$-module $\Val_k$ decomposes as  
\begin{displaymath}
 \Val_k \cong \bigoplus_\lambda \Gamma_\lambda,
\end{displaymath}
where $\lambda$ ranges over all highest weights such that $|\lambda_2| \leq 2$, $|\lambda_i| \neq 1$ for all $i$ and
$\lambda_i=0$ for $i > \min\{k,n-k\}$. In particular, it is multiplicity-free.
\end{Proposition}

\section{The Laplacian on the Grassmann manifold}

In this section $\pi:\SO(8)\rightarrow \Gr_k$ denotes the projection mapping each matrix  to the plane
spanned by its first $k$ columns. We also let $\mathbb S^1$ be the unit circle and define $\Phi:(\mathbb S^1)^4 \to
\SO(8)$ by \begin{displaymath}
\Phi(\theta_1,\ldots, \theta_4):=\left(\begin{matrix}
C&-S\\S&C
\end{matrix}\right) \in \SO(8),
\end{displaymath}
where
\begin{displaymath}
 C:=\left(\begin{matrix}
\cos\theta_1 & 0 & 0 & 0\\
0 & \cos\theta_2 & 0 & 0\\
0 & 0 & \cos\theta_3 &0\\
0 & 0 & 0 & \cos\theta_4
\end{matrix}\right),
\quad  S:=\left(\begin{matrix}
\sin\theta_1 & 0 & 0 & 0\\
0 & \sin\theta_2 & 0 & 0\\
0 & 0 &\sin\theta_3 & 0\\
0 & 0 & 0 & \sin\theta_4
\end{matrix}\right).
\end{displaymath} 

The image of $\Phi$ is a maximal torus of $\SO(8)$. We denote by $T$ the projection of this torus to $\Gr_k$,
which is a flat totally geodesic submanifold  of dimension $k$. By Corollaries 
\ref{coro_angles_3} and
\ref{coro_angles_4}, each $\spsp$-orbit has non-empty intersection with $T$. 

\begin{Proposition}
Each $\spsp$-orbit intersects $T$ orthogonally along a curve of the form
$c(t)=\pi\circ\Phi(\theta_1+t,\ldots,\theta_4+t)$; i.e. the tangent space to $T$ at $c(t)$ is spanned by $c'(t)$ and a
collection of vectors orthogonal to the orbit $\spsp\cdot c(t)$.
\end{Proposition}

\proof
By Corollary \ref{cor_param_angles}, the curve $c$ is contained in a single orbit. It remains to show that the
intersection of an orbit with $T$ is orthogonal.

Let us take the following basis of  $\mathfrak g=T_e\spsp$, viewed as a subspace of $\mathfrak{so}_8$:
\begin{equation}\label{basis}
\left(\begin{matrix}
0 & -\Id\\
\Id & 0
\end{matrix}\right),
\left(\begin{matrix} 
L_q & 0\\
0 & 0
\end{matrix}\right), 
\left(\begin{matrix} 
0 & 0\\
0 & L_q
\end{matrix}\right),
\left(\begin{matrix}
0 & L_q\\
L_q &0 
\end{matrix}\right),
\left(\begin{matrix} 
R_q & 0\\
0 & R_q
\end{matrix}\right),\qquad q=\mathbf i,\mathbf j,\mathbf k
\end{equation}
where $L_q,R_q\in \mathrm{End}_\R(\mathbb H)=\mathrm{End}_\R(\R^4)$ correspond to left and right multiplication by $q$
respectively. Let $N_i=\frac{\partial\Phi}{\partial \theta_i}-\frac{\partial\Phi}{\partial \theta_{i+1}}, 1 \leq
i \leq 3$, be bi-invariant vector fields defined on the maximal torus of $\SO(8)$. These vectors, together with the
vector $\sum_i \frac{\partial\Phi}{\partial \theta_i}$, span
the tangent space at each point of the maximal torus.

It is straightforward to check that $(N_i)_e$ is orthogonal to $\mathfrak g$, with respect to the Killing form of
$\mathfrak{so}_8$.  By
right-invariance, $(N_i)_g \bot \mathfrak g\cdot g$ for every $g$ in the maximal torus. Since $N_i \bot \ker d\pi$,
and $\pi$ is a riemannian submersion, we deduce that $(d\pi)_g N_i$ is orthogonal to the orbit $\spsp\cdot \pi(g)$.
 Since these vectors, together with $c'(t)$, span the tangent space of $T$ at $\pi(g)$, the statement follows. 
\endproof

Let $\vol:T \to \R$ be the function which assigns to $t \in T$ the volume of the orbit $\spsp\cdot
t$. By \cite[Corollary 1 and Proposition 1]{pacini03}, this function is positive and smooth on a dense subset of $T$.

\begin{Proposition}\label{prop_laplacian_level_sets}
Let $f$ be a smooth function on $\Gr_k$ which is invariant under $\spsp$. Let $\Delta$ be the
{Laplace-Beltrami
operator} acting on
smooth functions on $\Gr_k$. Let $\Delta_T$ be the Laplacian acting on functions on  $T$.
 Then, at all points where $\vol$ is strictly positive,
 \begin{displaymath}
  (\Delta f)|_T =\Delta_T f|_T - \langle \nabla (f|_T),\nabla(\log \vol)\rangle.
 \end{displaymath}
\end{Proposition}

\proof
By the previous proposition, there exists an orthonormal moving frame $E_1,\ldots, E_N$ on $\Gr_k$ such
that
$E_1,\ldots, E_d$ are orthogonal to the $\spsp$ orbits, and $E_1,\ldots, E_{k-1}$ span the tangent spaces of $T$.
Since $T$ is flat, we can assume that $\nabla_{E_i}E_j|_T=0$ for $i,j=1,\ldots,k$. Since $f$ is constant on the
orbits, 
\begin{displaymath}
 \nabla f = \sum_{i=1}^{k-1} E_i(f) E_i.
\end{displaymath}

Hence, on $T$,
\begin{align*}
 \Delta(f) &= -\mathrm{div}(\nabla f)\\
& =-\sum_j \sum_{i=1}^{k-1} \langle E_j,\nabla_{E_j}( E_i(f) E_i)\rangle\\
& =-\sum_{i=1}^{k-1} E_i\circ E_i(f)+\sum_{i=1}^{k-1} E_i(f)\sum_{ j=k}^{N} \langle
\nabla_{E_j}E_j,E_i\rangle\\
& =\Delta_T f+\langle \nabla f,\vec H\rangle,
\end{align*}
where $\vec H$ denotes the mean curvature vector of the $\spsp$-orbits. The result follows from the identity (cf. e.g.
\cite{pacini03})
\begin{displaymath}
 \vec H=-\nabla \log \vol.
\end{displaymath}
\endproof

\begin{Proposition}\label{prop_volume}
Let $g=\Phi(\theta_1,\ldots, \theta_k)$. The orbit $\spsp\cdot\pi(g)\subset \Gr_k$ has volume
\begin{align*}
\vol & ={c_2}\ \left|\sin(\theta_1-\theta_2)\right|^3 \cos(\theta_1-\theta_2)^2 & \mbox{ if } k=2,\\
\vol & ={c_3} \prod_{1\leq i<j\leq 3} \left| \sin(\theta_i-\theta_j)\right| \prod_{m\in \mathbb Z_3}
\left|\sin(\theta_{m+1}+\theta_{m+2}-2\theta_m)\right| & \mbox{ if } k=3,\\
\vol & ={c_4} \prod_{1\leq i<j\leq 4} \left| \sin(\theta_i-\theta_j)\right| \prod_{\{h,l\},\{m,n\}}
\left| \sin(\theta_h+\theta_l-\theta_m-\theta_n)\right| & \mbox{ if } k=4,
\end{align*}
where the last product runs over all unordered partitions $\{h,l\},\{m,n\}$ of $\{1,2,3,4\}$ into two disjoint pairs,
and $c_k$ is a constant depending only on $k$.
\end{Proposition}

\proof
We sketch the computation for $k=4$, the cases $k=2,3$ being similar. We just need to find the jacobian of the natural
map $\psi:\spsp\to \spsp\cdot\pi(g)$. By left-invariance, it is enough to compute $\mathrm{jac}(\psi)$ at $\mathfrak
g=T_e\spsp$. We
will use again the basis \eqref{basis} of $\mathfrak g$. The tangent space at $\pi(g)$ of $\Gr_4$ is
identified
using $d\pi\circ g^t$ with the horizontal part $\mathfrak m$ of $\mathfrak{so}_8$. This
way, for $X\in\mathfrak g$
\begin{displaymath}
 d\psi(X)=\pi_\mathfrak m(g^t X g)
\end{displaymath}
where $\pi_\mathfrak m\colon \mathfrak{so}_8\rightarrow \mathfrak m\equiv M_{4\times 4}(\R)$ consists of taking the
lower left block of the matrix. After identifying $\mathfrak m$ with $\R^{16}$, the matrix $A\in M_{13\times 16}(\R)$
associated with $d\psi$ is easily computed. The jacobian of $\psi$ is 
(up to constants) the determinant of $A$, with three rows of zeros removed. By suitably reordering the rows of $A$, one gets a structure of $4\times 4$ diagonal blocks, which makes the computation of the determinant an elementary task. 
\endproof

\begin{Proposition} \label{prop_tasaki_like_functions}
Let $f_{k,i}$ be the $\spsp$-invariant functions on $\Gr_k$ defined in the introduction. Then
\begin{align*}
\Delta(f_{k,0}) & = 0, \quad k=0,\ldots,4\\
\Delta(f_{2,1}) &= 28 f_{2,1}-12,\\
\Delta(f_{3,1}) &= 28 f_{3,1}-36,\\
\Delta(f_{3,2}) &= 60 f_{3,2}-34 f_{3,1}+18,\\
\Delta(f_{4,1}) &= 28 f_{4,1}-72,\\ 
\Delta(f_{4,2}) &= 40 f_{4,2}-2f_{4,1}-12,\\  
\Delta(f_{4,3}) &= 60 f_{4,3}+8 f_{4,2}-68 f_{4,1}+48,\\
\Delta(f_{4,4}) &= 96f_{4,4}+64 f_{4,1}-92 f_{4,3}-152 f_{4,2}+24.
\end{align*}
\end{Proposition}
\proof
It is enough to prove the identities on $T$.  By continuity, it suffices to prove them on the dense subset of
points corresponding to orbits of strictly positive volume. By Propositions \ref{prop_laplacian_level_sets}
and
\ref{prop_volume}, and using $\lambda_{ij}=\cos(\theta_i-\theta_j)$, this is a straightforward but lengthy computation. For instance,  $\Delta f_{2,1}$ is  computed by means of
\begin{displaymath}
 \Delta_T f_{2,1}=-4+8\cos^2(\theta_2-\theta_1),
\end{displaymath}
\begin{displaymath}
 \nabla f_{2,1}=2\cos(\theta_2-\theta_1)\sin(\theta_2-\theta_1)\left(\frac{\partial}{\partial \theta_1} -\frac{\partial}{\partial \theta_2}\right)
\end{displaymath}
\begin{displaymath}
 \nabla\log \vol=\frac{5\cos^2(\theta_2-\theta_1)-2}{\cos(\theta_2-\theta_1)\sin(\theta_2-\theta_1)}\left(-\frac{\partial}{\partial \theta_1} +\frac{\partial}{\partial \theta_2}\right).
\end{displaymath}

\endproof

\begin{Corollary} \label{cor_eigenvectors}
In each $\tilde \Gamma_\lambda$, there exists a unique (up to scale) invariant eigenfunction of the
Laplace-Beltrami operator on $\Gr_k$:
\begin{displaymath}
\renewcommand{\arraystretch}{1.3} 
 \begin{array}{c | c | c | c} 
k& \mbox{eigenfunction}& \mbox{eigenvalue} & \tilde \Gamma_\lambda \\  \hline \hline
0 & f_{0,0} & 0 & (0,0,0,0)\\ \hline
1 & f_{1,0} & 0 & (0,0,0,0)\\ \hline 
2 & f_{2,0} & 0 & (0,0,0,0)\\
2&7f_{2,1}-3f_{2,0}& 28 & (2,2,0,0)\\ \hline
3 & f_{3,0} & 0 & (0,0,0,0)\\
3 & 7f_{3,1}-9f_{3,0}  & 28 & (2,2,0,0)\\
3& 16 f_{3,2}-17f_{3,1}+15f_{3,0} & 60 & (4,2,2,0)\\ \hline
4 & f_{4,0} & 0 & (0,0,0,0)\\
4& 7f_{4,1}-18f_{4,0} & 28 & (2,2,0,0)\\
 4& 6f_{4,2}-f_{4,1} & 40 & (2,2,2,2)\\
4& 20f_{4,3}+8f_{4,2}-{43}f_{4,1}+66 f_{4,0} & 60 & (4,2,2,0)\\
4& 63 f_{4,4}-161 f_{4,3} -194 f_{4,2}+226 f_{4,1}-210f_{4,0}
&96 & (6,2,2,2)
 \end{array}
\renewcommand{\arraystretch}{1} 
\end{displaymath}
\end{Corollary}

\proof
To check that these functions are eigenvectors of the Laplacian with the given eigenvalues is easy using the previous
proposition. 

Let us show that these functions belong to $\tilde \Gamma_\lambda$ as stated in the last
column. 

It follows from Proposition \ref{prop_eigenvalues} that the eigenspaces corresponding to the eigenvalues $28$ and $60$
are given by $\tilde \Gamma_{(2,2,0,0)}$ and $\tilde \Gamma_{(4,2,2,0)}$. 

The eigenspace corresponding to the eigenvalue $40$ is given by $\tilde \Gamma_{(2,2,2,2)} \oplus \tilde
\Gamma_{(4,0,0,0)}$. The irreducible representation $\tilde \Gamma_{(4,0,0,0)}$ does not contain any
$\spsp$-invariant vector (otherwise $\dim \Val_1^{\spsp}$ would be larger than $1$, e.g. by Proposition
\ref{prop_val_decomposition}). Therefore an invariant 
eigenvector corresponding to the eigenvalue $40$ must belong to $\tilde \Gamma_{(2,2,2,2)}$. 

The eigenspace corresponding to the eigenvalue $96$ is given by $\tilde \Gamma_{(6,2,2,2)} \oplus \tilde
\Gamma_{(4,4,4,0)}$. The representation $\tilde \Gamma_{(4,4,4,0)}$ does not contain any $\spsp$-invariant vector. This
can
be checked using Weyl's character formula or a computer algebra system like LiE \cite{lie}. An
invariant eigenvector corresponding to the eigenvalue $96$ must thus belong to $\tilde \Gamma_{(6,2,2,2)}$.

Finally, to see that each $\tilde \Gamma_\lambda$ contains only one invariant function on $\Gr_k$, it is enough to
remark that each such function is the Klain function of an invariant valuation by Proposition
\ref{prop_val_decomposition}. By comparing dimensions (see
table \eqref{table_dims}), the claim follows. 
\endproof

Theorem \ref{thm_hadwiger_type} follows from Corollary \ref{cor_eigenvectors} and
Proposition \ref{prop_val_decomposition}. More precisely, each $\SO(8)$-representation $\tilde \Gamma_\lambda$ from the last column of the table enters the decomposition of $\Val_k$ by Proposition \ref{prop_val_decomposition}. By Schur's lemma and the injectivity of the Klain embedding, $\Val_k$ contains an $\spsp$-invariant valuation with the Klain function given in the second column.  Since these functions are linearly independent, we deduce from the dimensions in Table \ref{table_dims} that these valuations form a basis of $\Val_k^{\spsp}$.

Since we want to construct these valuations as explicitly as possible, we follow however a different path which allows
to compute Crofton measures
associated to the constructed valuations.

\section{Multipliers of the cosine transform}

Let $V \cong \R^n$ be a euclidean vector space. Set $\rho:=\frac{n}{2}$. The
$\alpha$-cosine transform $T_{k,k}^\alpha$ is defined for $\alpha \in \mathbb{C}$ with $\re
\alpha>\rho$ by  
\begin{align*}
 L^2(\Gr_k(\R^n)) & \to L^2(\Gr_k(\R^n))\\
 f & \mapsto \left[E \mapsto \int_{\Gr_k} f(F) |\cos(E,F)|^{\alpha-\rho} dF\right]
\end{align*}
and by meromorphic continuation for all $\alpha \in \mathbb{C}$.  
 
The case $\alpha=\rho+1$ yields the classical {\it cosine transform} \cite{lutwak90}, also denoted by $T_{k,k}$. 

Since $T_{k,k}^\alpha$ intertwines the $\SO(n)$-action, it acts as a scalar on each irreducible representation of
$\SO(n)$ which
enters the decomposition of  $L^2(\Gr_k(\R^n))$. The precise value of this constant was computed by \'Olafsson
and
Pasquale \cite{olafsson_pasquale} (compare also \cite{olaffson_pasquale_rubin} and \cite{zhang_genkai}). 

Let 
\begin{displaymath}
 \Gamma_k(\lambda):=\prod_{j=1}^k \Gamma\left(\lambda_j-\frac{j-1}{2}\right), \quad
\lambda=(\lambda_1,\ldots,\lambda_k) \in \mathbb{C}^k
\end{displaymath}
be the Siegel $\Gamma$-function. 

\begin{Theorem}[\'Olafsson-Pasquale]
 Let $\lambda=(\lambda_1,\ldots,\lambda_k)$ be a highest weight for $\SO(n)$ such that $\Gamma_\lambda$
enters the
decomposition of $L^2(\Gr_k(\R^n))$. Then $T_{k,k}^\alpha$ acts on $\Gamma_\lambda$ by
the scalar 
\begin{displaymath}
 c_{n,k}^\alpha:=(-1)^{\frac{|\lambda|}{2}} \frac{\Gamma_k(\rho)
\Gamma_k\left(\frac{\alpha-\rho+k}{2}\right)\Gamma_k\left(\frac{
-\alpha+\rho+\lambda}{2}\right)}{\Gamma_k\left(\frac{k}{2}\right)\Gamma_k\left(\frac{-\alpha+\rho}{2}
\right)\Gamma_k\left(\frac{\alpha+\rho+\lambda}{2}\right)}.
\end{displaymath}
In this formula, a complex number $z$ is identified with the vector $(z,\ldots,z) \in \mathbb{C}^k$. 
\end{Theorem}

\begin{Corollary}
 Let $\lambda=(\lambda_1,\ldots,\lambda_k,0,\ldots,0)$ be a highest weight of $\SO(n)$ such that $\Gamma_\lambda$
enters the decomposition of $\Val_k$ with $1 \leq k \leq \frac{n}{2}$. Then $T_{k,k}$ acts on
$\Gamma_\lambda$ by the scalar 
\begin{displaymath}
 c_{n,k}:=(-1)^{\frac{a}{2}-1} 
\frac{b'!(n-b'+1)!\Gamma\left(\frac{k+1}{2}\right)\Gamma\left(\frac{n-k+1}{2}\right)\Gamma\left(\frac{a-1}{2}\right)}{
2\pi n! \Gamma\left(\frac{n+1+a}{2}\right)}.
\end{displaymath}
Here $a:=\lambda_1$, $b$ is the depth of $\lambda$ (i.e. $\lambda_b \neq 0, \lambda_{b+1}=0$), and $b':=\max\{1,b\}$. 
\end{Corollary}

\proof
Clearly $\Gamma_k(\alpha)$ is well-defined and non-zero for $\alpha \in \R, \alpha > \frac{k-1}{2}$. We thus
have 
\begin{align*}
 c_{n,k} & =\lim_{\alpha \to \rho+1} c_{n,k}^\alpha\\
& = (-1)^{\frac{|\lambda|}{2}} \frac{\Gamma_k(\rho)\Gamma_k\left(\frac{k+1}{2}\right)}{ \Gamma_k\left(\frac{k}{2}\right)
\Gamma_k\left(\frac{n+1+\lambda}{2}\right)}
\lim_{\alpha \to \rho+1} \frac{
\Gamma_k\left(\frac{-\alpha+\rho+\lambda}{2}\right)}{\Gamma_k\left(\frac{-\alpha+\rho}{2}\right)} .
\end{align*}

Recall that, if $n$ is odd, we have $\lambda_j \in \{0,2\}$ for all $j>1$. If $n$ is even, then $\lambda_j \in \{0,2\}$ for $1<j<\frac{n}{2}$ and $\lambda_\frac{n}{2} \in \{0,2,-2\}$. 

Let us consider the first factor. Clearly 
\begin{displaymath}
  \frac{\Gamma_k\left(\frac{k+1}{2}\right)}{
\Gamma_k\left(\frac{k}{2}\right)}=\frac{\Gamma\left(\frac{k+1}{2}\right)}{\Gamma\left(\frac12\right)}.
\end{displaymath}

Next, we compute 
\begin{align*}
 \frac{\Gamma_k(\rho)}{\Gamma_k\left(\frac{n+1+\lambda}{2}\right)} & =
\frac{\Gamma\left(\frac{n-k+1}{2}\right)}{\Gamma\left(\frac{n+1+a}{2}\right)}
\prod_{j=2}^k \frac{\Gamma\left(\frac{n-j+2}{2}\right)}{\Gamma\left(\frac{n-j+2+\lambda_j}{2}\right)}.
\end{align*}

If $\lambda_j=0$, then the corresponding factor in the product equals $1$, while it equals $\frac{2}{n-j+2}$ if
$\lambda_j=2$. If $n$ is odd or $\lambda_\frac{n}{2} \neq -2$, the product thus equals $\frac{2^{b'-1}(n-b'+1)!}{n!}$. 

The last factor may be rewritten as
\begin{align*}
 \lim_{\alpha \to \rho+1} \frac{ \Gamma_k\left(\frac{-\alpha+\rho+\lambda}{2}\right)}{\Gamma_k\left(\frac{-\alpha+\rho}{2}\right)} & = \frac{\Gamma\left(\frac{a-1}{2}\right)}{\Gamma\left(-\frac12\right)}
\prod_{j=2}^k \lim_{x \to 0}
\frac{\Gamma\left(\frac{x+\lambda_j-j}{2}\right)}{\Gamma\left(\frac{x-j}{2}\right)}.
\end{align*}

If $\lambda_j=0$, then the corresponding term is $1$. If
$\lambda_j=2$, then the corresponding term equals
\begin{displaymath}
\lim_{x \to 0} \frac{\Gamma\left(\frac{x+2-j}{2}\right)}{\Gamma\left(\frac{x-j}{2}\right)}  = - \frac{j}{2}. 
\end{displaymath}

If $\lambda_\frac{n}{2} \neq -2$, we thus get that 
\begin{displaymath}
\lim_{\alpha \to \rho+1} \frac{ \Gamma_k\left(\frac{-\alpha+\rho+\lambda}{2}\right)}{\Gamma_k\left(\frac{-\alpha+\rho}{2}\right)} = \frac{\Gamma\left(\frac{a-1}{2}\right)}{\Gamma\left(-\frac12\right)}
\frac{(-1)^{b'-1}
b'!}{2^{b'-1}}=\frac{\Gamma\left(\frac{a-1}{2}\right) b'! (-1)^{b'}}{\sqrt{\pi} 2^{b'}}.
\end{displaymath}

Putting these pieces together yields for $\lambda_\frac{n}{2} \neq -2$ 
\begin{displaymath}
 c_{n,k}=(-1)^{\frac{a}{2}-1} 
\frac{b'!(n-b'+1)!\Gamma\left(\frac{k+1}{2}\right)\Gamma\left(\frac{n-k+1}{2}\right)\Gamma\left(\frac{a-1}{2}\right)}{
2\pi n! \Gamma\left(\frac{n+1+a}{2}\right)}.
\end{displaymath}

Finally, if $n$ is even, let us compare the cases $(a,2,\ldots,2,2)$ and $(a,2,\ldots,2,-2)$. The first factor gets
multiplied by $\frac{\Gamma\left(\frac{n}{4}+2\right)}{\Gamma\left(\frac{n}{4}\right)}$, while the second factor gets
multiplied by  $\frac{\Gamma\left(\frac{n}{4}\right)}{\Gamma\left(\frac{n}{4}+2\right)}$. Hence the constant $c_{n,k}$
is the same in both cases, which completes the proof. 
\endproof

\begin{Corollary}
The cosine transform acts by the following scalars 
\begin{displaymath}
\renewcommand{\arraystretch}{1.3} 
\begin{array}{c | c | c} 
k & \tilde \Gamma_\lambda & c\\ \hline 
2 & (0,0,0,0) & \frac{1}{7}\\
2 & (2,2,0,0) & \frac{1}{252}\\ \hline
3 & (0,0,0,0) & \frac{32}{105 \pi}\\
3 & (2,2,0,0) &  \frac{8}{945 \pi} \\
3 & (4,2,2,0) & -\frac{8}{24255\pi}\\ \hline
4 & (0,0,0,0) & \frac{3}{35}\\
4 & (2,2,0,0) &  \frac{1}{420}\\
4 & (2,2,2, 2) & \frac{1}{1470}\\
4 & (4,2,2,0) & -\frac{1}{10780}\\
4 & (6,2,2, 2) & \frac{1}{70070} 
\end{array} 
\renewcommand{\arraystretch}{1} 
\end{displaymath}
\end{Corollary}

\section{Construction of invariant valuations}

\begin{Proposition} \label{prop_hadwiger_eigenvectors}
 There exist valuations in $\Val_k^{\spsp}$, $k=0,\ldots,8$, whose Klain functions on $\Gr_k
\cong \Gr_{\min\{k,8-k\}}$ are given by the eigenfunctions from Corollary
\ref{cor_eigenvectors}. These valuations form a basis of  $\Val_k^{\spsp}$.
\end{Proposition}

\proof
Let $g \in C(\Gr_k)$ and define a valuation in $\mu \in \Val_k^+$ by 
\begin{displaymath}
\mu(K):=\int_{\Gr_k} g(E) \vol(\pi_EK) dE,
\end{displaymath}
{where $\pi_E:\h^2 \to E$ is the orthogonal projection. }
Then $\Kl_\mu = T_{k,k}g$. 

If $f$ is an eigenfunction from the table in Corollary \ref{cor_eigenvectors}, then the cosine transform $T_{k,k}$ acts
by a non-zero scalar $c$. Setting $g:=c^{-1} f$ we get $\Kl_\mu=f$. 

By looking at their Klain functions, we deduce that the so-constructed valuations are linearly independent in each
degree of homogeneity. By comparing with the dimensions in  \eqref{table_dims}, they actually must form a basis. 
\endproof

\proof[Proof of Theorem \ref{thm_hadwiger_type}]
The theorem follows from Proposition \ref{prop_hadwiger_eigenvectors} by noting that the transformation matrix between
the $f_{k,i}$ and the eigenvectors is invertible. 
\endproof


\def\cprime{$'$}

\end{document}